\documentclass[11pt,leqno,a4paper]{article}

\usepackage{amsmath, amssymb, amsthm, amsfonts}
\usepackage{mathrsfs}
\usepackage{color}
\usepackage{hyperref}
\theoremstyle{plain}
\newtheorem{thm}{Theorem}

\newtheorem{lem}{Lemma}

\theoremstyle{remark}
\newtheorem{rem}{Remark}

\renewcommand{\Re}{{\rm Re\,}}

\def\K{\mathop{\mbox{\bf\Large K}}}

\numberwithin{equation}{section}

\begin{document}

\title{Multiple-correction and summation of the rational series}
\author{Xiaodong Cao*, Cristinel Mortici}
\date{}

\maketitle
\footnote[0]{* Corresponding author.}

\footnote[0]{\emph{E-mail address}: caoxiaodong@bipt.edu.cn~(X.D. Cao), cristinel.mortici@hotmail.com~(C. Mortici) }


\footnote[0]{This work is supported by
the National Natural Science Foundation of China (Grant No.11171344) and the Natural
Science Foundation of Beijing (Grant No.1112010).}

\footnote[0]{
Xiaodong Cao:
Department of Mathematics and Physics,
Beijing Institute of Petro-Chemical Technology,
Beijing, 102617, P. R. China \\
Cristinel Mortici:
Dept. of Mathematics,
Valahia University of T\^argovi\c{s}te, Bd. Unirii 18, 130082 Targovi\c{s}te, Romania\\
Academy of Romanian Scientists, Splaiul Independen\c{t}ei nr. 54, 050094 Bucharest, Romania\\
E-mail addresses: cristinel.mortici@hotmail.com}

\begin{abstract} The goal of this work is to formulate  a systematical method for looking for the simple closed form or continued fraction representation of a class of rational series. As applications, we obtain the continued fraction representations for the alternating Mathieu series and some rational series. The main tools are multiple-correction and two of Ramanujan's continued fraction formulae involving the quotient of the gamma functions.
\end{abstract}

\section{Introduction}

Let the general term  of an infinite series have the form
\begin{align}
u_n=\frac{P_l(n)}{Q_m(n)},
\end{align}
where $P_l(n), Q_m(n)$ are polynomials of degree $l$ and $m$ with real coefficients, respectively. Finding the sum of a rational series $\sum u_n$ in \emph{simple closed form} is a very important research area, see, e.g.,
R.L. Graham, D.E. Knuth and O. Patashnik~\cite{GKP}, M. Petkov\u{s}ek, H.S. Wilf and D. Zeilberger~\cite{PWZ}, H.S. Wilf~\cite{Wilf} and references therein. Throughout the paper, the \emph{simple closed form} always means a rational function with real coefficients. Otherwise, we expect that a continued fraction representation for the series may be discovered, see Chapter 12 in Berndt~\cite{Ber}, L. Lorentzen and H. Waadeland~\cite{LW}, or A. Cuyt, V.B. Petersen, B. Verdonk, H. Waadeland, W.B. Jones~\cite{CPV}. The main purpose of this paper is to investigate a kind of fundamental rational series in a unified setting, which contains some mathematical constants and series, such as Catalan constant, $\zeta(2)$, Ap\'ery number $\zeta(3)$, the Mathieu series and the alternating Mathieu series, etc. The Mathieu series was introduced by \'Emile Leonard Mathieu in his book~\cite{Mathieu}, which is defined by
\begin{align}
S(r):=\sum_{m=1}^{\infty}\frac{2m}{(m^2+r^2)^2},\quad (r>0),\label{Mathieu-Series}
\end{align}
while the \emph{alternating Mathieu series} is given as follows
\begin{align}
\tilde{S}(r):=\sum_{m=1}^{\infty}(-1)^{m-1}\frac{2m}{(m^2+r^2)^2}.
\quad (r>0).\label{Alternating Mathieu-Series}
\end{align}
The Mathieu series has important applications in science, such as in the theory of elasticity of solid bodies~\cite{Eme}, or in the problem of the rectangular
plate~\cite{Sch} and it is closely related to the Riemann Zeta function $\zeta$~\cite{CS}. Moreover, A. Jakimovski and D.C. Russell~\cite{JR} showed that an
extended form of the Mathieu series plays a role in examining Mercerian theorems for Ces\'aro summability. For the research history of the Mathieu series and the related series, readers interested may refer to R. Frontczak~\cite{Fro},
G.V. Milovanovi\'c and T.K. Pog\'any~\cite{MP}, C. Mortici~\cite{Mor-M}, T.K. Pog\'any, H.M. Srivastava and \u{Z}. Tomovski~\cite{PST}, and references therein.

The paper is organized as follows. In Sec.~2, we gives several notations and definitions for later use. In Sec.~3, based on the work in \cite{CY}, we shall formulate a systematical method to look for either a \emph{simple closed form} solution or the fastest possible finite continued fraction approximation solution for a class of linear difference equation of order one, sometimes we may guess further its continued fraction solution. In order to prove our new ``conjectures", in Sec. 4 we shall prepare some important tools, which are two of the famous Ramanujan's continued fraction formulas for the quotient of the gamma functions. In Sec.~5 and 6, we shall investigate the rational series $\sum\frac{1}{Q_l(n)}$ for $l=2,3$, respectively, where $Q_l(x)$ is a polynomial of degree $l$ in $x$. In Sec.~7, we shall continue to study two extended Mathieu series. In Sec.~8, we shall give two applications of main results in Sec.~7. For example, we establish first a new representation for the alternating Mathieu series in the form of  a linear combination of two continued fractions. In the last section, we analyze the related perspective of research in this direction.
\bigskip

\section{Notation and definition}
Throughout the paper, we use the digamma notation $\psi(z)=\frac{\Gamma'(z)}{\Gamma(z)}$. The set $\mathbb{Z}\backslash\mathbb{N}$ means $\{0,-1,-2,\ldots\}$. The notation
$P_k(x)$~(or $Q_k(x)$) means a polynomial of degree $k$
in $x$, while $U(x)$~(or $V(x)$) denotes a rational function in $x$ . We shall use the $\Phi(k;x)$ to denote a polynomial of degree $k$ in $x$ with the leading coefficient equals one, which may be different at each occurrence. Let $(a_n)_{n\ge 1}$ and $(b_n)_{n\ge 0}$ be two sequences of real or complex numbers. The generalized continued fraction
\begin{align}
\tau=b_0+\frac{a_1}{b_1+\frac{a_2}{b_2+\ddots}}=b_0+
\begin{array}{ccccc}
a_1 && a_2 &       \\
\cline{1-1}\cline{3-3}\cline{5-5}
 b_1 & + & b_2 & + \cdots
\end{array}
=b_0+\K_{n=1}^{\infty}
\left(\frac{a_n}{b_n}\right)
\end{align}
is defined as the limit of the $n$th approximant
\begin{align}
\frac{A_n}{B_n}=b_0+\K_{k=1}^{n}\left(\frac{a_k}{b_k}\right)
\end{align}
as $n$ tends to infinity. In line with Ramanujan we adopt the convention that if $a_N=0$ and $a_n\neq 0$ for all $n<N$, then the
continued fraction $\tau$ \emph{terminates} and  has the value
\begin{align}
\tau=\frac{A_{N-1}}{B_{N-1}}=
b_0+\K_{k=1}^{N-1}\left(\frac{a_k}{b_k}\right).
\end{align}
This has the advantage that the continued fraction is always well defined. For the continued fraction theory, see L. Lorentzen and H. Waadeland~\cite{LW}, and A. Cuyt, V.B. Petersen, B. Verdonk, H. Waadeland, W.B. Jones~\cite{CPV}, or other classical books therein.
For the sake of convenience, the sum $\sum_{j=n_1}^{n_2}a_j$ or the continued fraction $\K_{j=n_1}^{n_2}\left(\frac{c_j}{d_j}\right)$ for $n_2<n_1$ is stipulated to be zero.

\bigskip

In order to describe our method clearly, Let us recall three definitions introduced in~\cite{CTZ}.
\bigskip

\noindent {\textbf {Definition 1.}}
Let $f(x)$ be a function defined on $(x_0,+\infty)$ for some real $x_0$. We assume that there exists a fixed positive number $\mu$ and a constant $c\neq 0$ such that $\lim_{x\rightarrow +\infty}x^{\mu}f(x)=c$. We define
\begin{align}
 \mathrm{R}(f(x)):=\mu,
\end{align}
where $\mu$ is the exponent of $x^{\mu}$. For convenience, $\mathrm{R}(0)$ is stipulated to be infinity. Hence, $\mathrm{R}(f(x))$ characterizes the rate of convergence for $f(x)$ as $x$ tends to infinity.

\bigskip

\noindent {\textbf {Definition 2.}} Let $c_0\neq 0$, and $x$ be a free variable. Let  $(a_n)_{n=1}^{\infty}$, $(b_n)_{n=1}^{\infty}$ and $(c_n)_{n=1}^{\infty}$ be three real sequences. The formal continued fraction
\begin{align}
\frac{c_0}{\Phi(\nu;x)+\K_{n=1}^{\infty}\left(\frac{a_n}{x+b_n}
\right)}
\end{align}
is said to be a \emph{Type-I} continued fraction, i.e. when $n\ge 2$ the $n$-th partial denominator is a linear function in $x$ . While,
\begin{align}
\frac{c_0}{\Phi(\nu;x)+\K_{n=1}^{\infty}\left(\frac{a_n}{x^2+b_n x+c_n}\right)}
\end{align}
is said to be a \emph{Type-II} continued fraction, i.e. when $n\ge 2$  the $n$-th partial denominator is a polynomial of degree 2 in $x$.
\bigskip

\begin{rem} The \emph{Type-I} and \emph{Type-II} are two kinds of fundamental structures we often meet. Similarly, we may define other type continued fractions. Certainly, there also exists the hybrid-type continued fractions. In this paper, we shall not discuss the involved problems. It should be remarked that for the formal power series solution of the equation ~\eqref{Difference Equation} below, its structure is unique. However, the structure of the formal continued fraction solution is more complicated, which may be a \emph{Type-I} or \emph{Type-II}, or a \emph{hybrid-type}. This is the main motivation that we introduce the definitions of a \emph{Type-I, Type-II}, or other type to classify the fastest possible continued fraction solution.
\end{rem}
\bigskip

\noindent {\textbf {Definition 3.}} If the sequence $(b_n)_{n=1}^{\infty}$ is a constant sequence $(b)_{n=1}^{\infty}$ in the \emph{Type-I}~( or \emph{Type-II})
continued fraction, we call the number $\omega=b$~(or $\omega=\frac b2$) the $\mathrm{MC}$-point for the corresponding continued fraction.
We use $\hat{x}=x+\omega$ to denote the $\mathrm{MC}$-shift of $x$.
\bigskip

If there exists the $\mathrm{MC}$-point, we have the following \emph{simplified form}
\begin{align}
\frac{c_0}{\Phi_1(\nu;\hat{x})+\K_{n=1}^{\infty}\left(\frac{a_n}
{\hat{x}}\right)}\quad
\mbox{or}\quad
\frac{c_0}{\Phi_1(\nu;\hat{x})+\K_{n=1}^{\infty}\left(\frac{a_n}
{\hat{x}^2+d_n}\right)},\label{canonical-form}
\end{align}
where $d_n=c_n-\frac{b^2}{4}.$

\section{The multiple-correction method and the continued fraction solution of linear difference equation of order one}
Let $U(x)$ and $V(x)$ be rational functions in $x$. We consider the following linear difference equation of order one
\begin{align}
y(x)-U(x)y(x+1)-V(x)=0.\label{Difference Equation}
\end{align}
We are concern with the continued fraction solution of the difference equation~\eqref{Difference Equation} when $x$ tends to infinity~(or for ``large" $x$). One of our purpose is to try to look for the fastest possible finite continued
fraction approximation solution or guess the formal continued fraction solution.
In the other words, we are looking for the solution~(or approximation solution) of the equation~\eqref{Difference Equation}
in the \emph{formal continued fraction space}, which contains the subspace
of the \emph{rational function fields}. Our method may be described as the following five steps.
\bigskip

\noindent \textbf{(Step 1)} Let us develop further the previous \emph{multiple-correction method} formulated in~\cite{CY}.
In fact, the
\emph{multiple-correction method} is a recursive algorithm, and one of its advantages is that by repeating correction-process we always can accelerate the rate of approximation. More precisely, every non-zero coefficient plays an important role in accelerating the rate of approximation. For the sake of completeness, we shall give a description in details. The \emph{multiple-correction method} consists of the following several steps.
\bigskip

\noindent {\bf(Step 1-1) The initial-correction.} The choice of initial-correction is
vital. Determine the initial-correction $\mathrm{MC}_0(x)=\frac{c_0}{\Phi_0(\nu;x)}$~(or $\mathrm{MC}_0(x)=c_0\Phi_0(\nu;x)$) such that
\begin{align}
&\mathrm{R}\left(\frac{c_0}{\Phi_0(\nu;x)}-U(x)
\frac{c_0}{\Phi_0(\nu;x+1)}-V(x)
\right)\\
=&
\max_{c, \Phi(\nu;x)}\mathrm{R}
\left(\frac{c}{\Phi(\nu;x)}-U(x)
\frac{c}{\Phi(\nu;x+1)}-V(x)\right).\nonumber
\end{align}
In the second case, we can modify the approach above as follows
\begin{align}
&\mathrm{R}\left(c_0~\Phi_0(\nu;x)-U(x)
c_0~\Phi_0(\nu;x+1)-V(x)
\right)\\
=&
\max_{c, \Phi(\nu;x)}\mathrm{R}
\left(c~\Phi(\nu;x)-U(x)c~\Phi(\nu;x+1)
-V(x)\right).\nonumber
\end{align}
In the sequel, we only give the exact expressions for the first case. For the second case, the correction-functions can be constructed mutatis mutandis as for the first case.

After determining the initial-correction, we define the  initial-correction error function $E_0(x)$ by
\begin{align}
E_0(x)=\mathrm{MC}_0(x)-U(x)\mathrm{MC}_0(x+1)-V(x).
\end{align}

Find $\mathrm{R}(E_0(x))$. If $E_0(x)\equiv 0$, then the difference equation~\eqref{Difference Equation} has a \emph{simple closed form} solution $\mathrm{MC}_0(x)$.
\bigskip

Now we explain how to determine all the coefficients in $\mathrm{MC}_0(x)$. Firstly, we try to look for $c_0$  and $\nu$, which satisfy the following condition
\begin{align}
\max_{c,l}\mathrm{R}
\left(\frac{c}{x^l}-U(x)\frac{c}{(x+1)^l}
-V(x)\right).
\end{align}
Secondly, If $\nu>0$, then we take $\mathrm{MC}_0(x)$ in the form $\frac{c_0}{\Phi_0(\nu;x)}$. Otherwise, we choose $\mathrm{MC}_0(x)=c_0 \Phi_0(-\nu;x)$. Thirdly, we may determine other coefficients in $\mathrm{MC}_0(x)$ by successively solving a linear equation.

\bigskip
\noindent {\bf(Step 1-2) The first-correction.} If there exists a real number $\kappa_1$ such that
\begin{align}
\mathrm{R}\left(\frac{c_0}{\Phi_0(\nu;x)+\frac{\kappa_1}{x}}
-U(x)
\frac{c_0}{\Phi_0(\nu;x+1)+\frac{\kappa_1}{x+1}}-V(x)
\right)
>\mathrm{R}\left(E_0(x)\right),
\end{align}
then we take the first-correction $\mathrm{MC}_1(x)=\frac{\kappa_1}{x+\lambda_1}$ with
\begin{align}
&\mathrm{R}\left(\frac{c_0}{\Phi_0(\nu;x)
+\frac{\kappa_1}{x+\lambda_1}}
-U(x)
\frac{c_0}{\Phi_0(\nu;x+1)+\frac{\kappa_1}{x+1+\lambda_1}}
-V(x)\right)\\
=&\max_{\lambda}\mathrm{R}\left(\frac{c_0}{\Phi_0(\nu;x)
+\frac{\kappa_1}{x+\lambda}}
-U(x)
\frac{c_0}{\Phi_0(\nu;x+1)+\frac{\kappa_1}{x+1+\lambda}}-V(x)\right).
\nonumber
\end{align}
Otherwise, we take the first-correction $\mathrm{MC}_1(x)=\frac{\kappa_1}{x^2+\lambda_{1,1}x
+\lambda_{1,2}}$ such that
\begin{align}
\mathrm{R}
&\left(\frac{c_0}{\Phi_0(\nu;x)
+\frac{\kappa_1}{x^2+\lambda_{1,1} x+\lambda_{1,2}}}
-U(x)
\frac{c_0}{\Phi_0(\nu;x+1)+\frac{\kappa_1}{(x+1)^2+\lambda_{1,1} (x+1)+\lambda_{1,2}}}-V(x)
\right)\\
=&
\max_{\kappa,\lambda_1,\lambda_2}\mathrm{R}
\left(\frac{c_0}{\Phi_0(\nu;x)
+\frac{\kappa}{x^2+\lambda_1 x+\lambda_2}}
-U(x)
\frac{c_0}{\Phi_0(\nu;x+1)+\frac{\kappa}{(x+1)^2+\lambda_1 (x+1)+\lambda_2}}-V(x)
\right).\nonumber
\end{align}
If $\kappa_1=0$, we need to stop the correction-process, which means that the rate of approximation can not be further improved only by making use of \emph{Type-I} or \emph{Type-II} continued fraction structure. In the other words, in order to improve the rate of approximation, we have to choose a more general continued fraction structure instead of it. More precisely, we take first $\mathrm{MC}_1(x)=\frac{\kappa_1}{x^j}$.
Then we need to begin from $j=1$ and try step by step. Once we have found that the convergence rate can be improved for the first positive integer, say $j_0$, we use $\Phi(j_0; n)$ to replace $n^{j_0}$ immediately, and determine all the corresponding coefficients of the polynomial $\Phi(j_0; n)$, which is the main new idea introduced in~\cite{CXY,Cao1}. Lastly, we choose $\mathrm{MC}_1(x)=\frac{\kappa_1}{\Phi(j_0; n)}$. In what follows,
we only describe our method for the structures of the \emph{Type-I}
and \emph{Type-II} continued fraction approximation.

Now we define the first-correction error function $E_1(x)$ by
\begin{align}
E_1(x)=\frac{c_0}{\Phi_0(\nu;x)
+\mathrm{MC}_1(x)}-U(x)\frac{c_0}{\Phi_0(\nu;x+1)
+\mathrm{MC}_1(x+1)}-V(x).
\end{align}
Find $\mathrm{R}(E_1(x))$. If $E_1(x)\equiv 0$, then the difference equation~\eqref{Difference Equation} has a \emph{simple closed form} solution
\begin{align}
\frac{c_0}{\Phi_0(\nu;x)
+\mathrm{MC}_1(x)}=
\begin{array}{cc}
c_0&  \\
\cline{1-1}
\Phi_0(\nu;x) & +
\end{array}\mathrm{MC}_1(x).
\end{align}

\bigskip

\noindent {\bf(Step 1-3) The second-correction to the $k$th-correction.} If $\mathrm{MC}_1(x)$ has the form
\emph{Type-I}, we take the second-correction
\begin{align}
\mathrm{MC}_2(x)=\frac{\kappa_1}
{x+\lambda_1+\frac{\kappa_2}{x+\lambda_2}}=
\begin{array}{ccc}
\kappa_1&  &\kappa_2\\
\cline{1-1}\cline{3-3}
x+\lambda_1 & +&x+\lambda_2
\end{array}
,
\end{align}
which satisfies
\begin{align}
\max_{\kappa,\lambda}\mathrm{R}
\left(\frac{c_0}{\Phi_0(\nu;x)
+\frac{\kappa_1}{x+\lambda_1+\frac{\kappa}{x+\lambda}}}
-U(x)
\frac{c_0}{\Phi_0(\nu;x+1)+\frac{\kappa_1}{x+1+\lambda_1
+\frac{\kappa}{x+1+\lambda}}}-V(x)
\right).
\end{align}
Similarly to the first-correction, if $\kappa_2=0$, we stop the correction-process.

If $\mathrm{MC}_1(x)$ has the form \emph{Type-II}, we take the second-correction
\begin{align}
\mathrm{MC}_2(x)=\frac{\kappa_1}
{x^2+\lambda_{1,1}x+\lambda_{1,2}+\frac{\kappa_2}
{x^2+\lambda_{2,1}x+\lambda_{2,2}}}
\end{align}
such that
\begin{align}
&\max_{\kappa,\lambda_{1},\lambda_{2}}\mathrm{R}
\left(\frac{c_0}{\Phi_0(\nu;x)
+\frac{\kappa_1}
{x^2+\lambda_{1,1}x+\lambda_{1,2}+\frac{\kappa}{x^2+\lambda_1 x+\lambda_2}}}\right.\\
&\quad\left.
-U(x)
\frac{c_0}{\Phi_0(\nu;x+1)+\frac{\kappa_1}
{(x+1)^2+\lambda_{1,1}(x+1)+\lambda_{1,2}+\frac{\kappa}{(x+1)^2
+\lambda_1 (x+1)+\lambda_2}}}-V(x)
\right).\nonumber
\end{align}
If $\kappa_2=0$, we also need to stop the correction-process.

Now we define the second-correction error function $E_2(x)$ by
\begin{align}
E_2(x)=\frac{c_0}{\Phi_0(\nu;x)
+\mathrm{MC}_2(x)}-U(x)\frac{c_0}{\Phi_0(\nu;x+1)
+\mathrm{MC}_2(x+1)}-V(x).
\end{align}
If $E_2(x)\equiv 0$, then the difference equation~\eqref{Difference Equation} has a \emph{simple closed form} solution
\begin{align}
\frac{c_0}{\Phi_0(\nu;x)
+\mathrm{MC}_2(x)}=
\begin{array}{cc}
c_0&  \\
\cline{1-1}
\Phi_0(\nu;x) & +
\end{array}\mathrm{MC}_2(x).
\end{align}
\bigskip

If we can continue the above correction-process to determine the $k$th-correction function $\mathrm{MC}_k(x)$ until some $k^*$ you want, then one may use a recurrence relation to determine the $k$th-correction $\mathrm{MC}_k(x)$. More precisely, in the case of \emph{Type-I} we choose
\begin{align}
\mathrm{MC}_k(x)=\K_{j=1}^{k}
\left(\frac{\kappa_j}{x+\lambda_j}\right)
=\mathrm{MC}_{k-1}(x)\begin{array}{cc}
& \kappa_k \\
\cline{2-2}
+ & x+\lambda_k
\end{array}
\end{align}
such that
\begin{align}
\max_{\kappa,\lambda}\mathrm{R}
\left(\frac{c_0}{\Phi_0(\nu;x)
+G(\kappa,\lambda;x)}
-U(x)
\frac{c_0}{\Phi_0(\nu;x+1)+G(\kappa,\lambda;x+1)}
-V(x)\right),
\end{align}
where
\begin{align*}
G(\kappa,\lambda;x):=\mathrm{MC}_{k-1}(x)\begin{array}{cc}
& \kappa \\
\cline{2-2}
+ & x+\lambda
\end{array}=
\begin{array}{ccccc}
\kappa_1 & & \kappa_{k-1} &      & \kappa \\
\cline{1-1}\cline{3-3}\cline{5-5}
x+\lambda_1 & +\cdots+ & x+\lambda_{k-1} & + &x+\lambda
\end{array}.
\end{align*}
While, in the case of \emph{Type-II} we take
\begin{align}
\mathrm{MC}_k(x)=\K_{j=1}^{k}\left(\frac{\kappa_j}
{x^2+\lambda_{j,1}x+\lambda_{j,2}}\right),
\end{align}
which satisfies
\begin{align}
\max_{\kappa,\lambda_{1},\lambda_{2}}\mathrm{R}
\left(\frac{c_0}{\Phi_0(\nu;x)
+H(\kappa,\lambda_1,\lambda_2;x)}
-U(x)
\frac{c_0}{\Phi_0(\nu;x+1)+H(\kappa,\lambda_1,\lambda_2;x+1)}
-V(x)\right),
\end{align}
where
\begin{align*}
H(\kappa,\lambda_1,\lambda_2;x):=&
\mathrm{MC}_{k-1}(x)\begin{array}{cc}
& \kappa \\
\cline{2-2}
+ & x^2+\lambda_1 x+\lambda_2
\end{array}\\
=&\begin{array}{ccccc}
\kappa_1 & & \kappa_{k-1} &      & \kappa \\
\cline{1-1}\cline{3-3}\cline{5-5}
x^2+\lambda_{1,1}x+\lambda_{1,2} & +\cdots+ & x^2+\lambda_{k-1,1}x+\lambda_{k-1,2} & + &x^2+\lambda_1 x+\lambda_2
\end{array}.
\end{align*}
In both cases, if $\kappa_k=0$, we have to stop the correction-process.

\bigskip

Now we define the $k$th-correction error function $E_k(x)$ by
\begin{align}
E_k(x)=\frac{c_0}{\Phi_0(\nu;x)
+\mathrm{MC}_k(x)}-U(x)\frac{c_0}{\Phi_0(\nu;x+1)
+\mathrm{MC}_k(x+1)}-V(x).\label{Error-term-Def}
\end{align}

Find $\mathrm{R}(E_k(x))$. Lastly, if $E_k(x)\equiv 0$, then the difference equation~\eqref{Difference Equation} has a \emph{simple closed form} solution
\begin{align}
\frac{c_0}{\Phi_0(\nu;x)
+\mathrm{MC}_k(x)}.
\end{align}

\bigskip

For the reader's convenience, we would like to give the complete \textbf{Mathematica} program for determining
all the coefficients in  $\mathrm{MC}_k(x)$.
\bigskip

{\bf (1)}. First, let the function $Ek[x]$ be defined by~\eqref{Error-term-Def}.
\bigskip

{\bf (2)}. Then we manipulate the following \textbf{Mathematica} command to
expand $Ek[x]$ into a power series in terms of $1/x$:
\begin{align}
\text{Normal}[\text{Series}
[Ek[x]
\text{/.}~ x\rightarrow
1/u, \{u,0,l_k\}]]\text{/.}~ u\rightarrow 1/x~(\text{// Simplify}).
\label{Ek-Mathematica-Program}
\end{align}
We remark that the variable $l_k$ needs to be suitable chosen according to the different functions and $k$. Another approach is that putting the whole thing over a common denominator such that $\mathrm{R}\left(E_k(x)\right)$
is strictly decreasing function of $k$. We may manipulate \textbf{Mathematica} commands ``Together" and ``Collect" to achieve them.
\bigskip

{\bf (3)}. Taking out the first some coefficients in the above power series, then we enforce them to be zero, and finally solve the related coefficients successively.

\begin{rem}
Actually, once we have found $\mathrm{MC}_k(x)$, \eqref{Ek-Mathematica-Program} can be used again to determine $\mathrm{R}\left(E_k(x)\right)$. In addition, we can apply it to check the general term formula for $\mathrm{MC}_k(x)$.
\bigskip

A lot of experiment results show that  $\mathrm{R}\left(E_k(x)\right)\ge\mathrm{R}\left(E_0(x)\right)+2k$ in the case of \emph{Type-I}, and $\mathrm{R}\left(E_k(x)\right)\ge\mathrm{R}\left(E_0(x)\right)+4k$ in the case of \emph{Type-II}, respectively. We also note that the $k$th-approximation solution
$
\frac{c_0}{\Phi_0(\nu;x)
+\mathrm{MC}_k(x)}
$
of the equation~\eqref{Difference Equation} may be written in the form $\frac{P(l_1;x)}{Q(l_2;x)}$ with $l_2=\nu+k$ or $l_2=\nu+2k$, respectively, which explains that our method provides indeed an effective approach for the approximation solution problem of the equation ~\eqref{Difference Equation}. We think that it should be the best possible rational approximation solution.
\bigskip
\end{rem}
\bigskip

\noindent \textbf{(Step 2)} Find the general term formula of the $k$th-correction. Here we often use some tools in number theory, difference equation, etc.

If one can not find the general term formulas of both the $n$-th partial numerator and denominator, then only the finite continued fraction approximation solution can be provided. For instance, the \emph{BBP-type} series of some mathematical constants like $\pi$, Catalan constant, $\pi^2$, etc.~(e.g. see \cite{CY}). At the same times, it predicts that finding a continued fraction representation started started from this series seems ``hopeless", one should replace other series expressions to try again and again, or find a linear combination solution of several continued fractions. Perhaps, an unexpected surprise will happen~!

\bigskip

\noindent \textbf{(Step 3)} If we are lucky, we find that $E_{k^*}(x)\equiv 0$ for some integer $k^*$, then we attain a \emph{simple closed form}~(or a finite continued fraction) solution of the difference equation~\eqref{Difference Equation}. Now we give two examples to illustrate our method.

\bigskip

\noindent{\textbf{Example 1}} Consider the following equation
\begin{align}
y(x)-y(x+1)=\frac{12 x^4-1}{(4 x^4+1  )^2}.\label{Example-1}
\end{align}
As $\mathrm{R}\left(\frac{12 x^4-1}{(4 x^4+1  )^2}\right)=4$, firstly, we look for the exponent $\nu$ and the $c_0$ such that $\nu, c_0$ satisfy the following condition
\begin{align*}
\max_{\nu,c}\mathrm{R}\left(\frac{c}{x^{\nu}}
-\frac{c}{(x+1)^{\nu}}-\frac{12 x^4-1}{(4 x^4+1  )^2}\right).
\end{align*}
It is not difficult to see that $\nu=3$. Then
we put three functions in the expression above over a common denominator, then let the first coefficient $-12 + 48 c$ in the numerator to be zero. In this way, we find $c_0=\frac 14$. Now we
choose the first-correction $\mathrm{MC}_0(x)=\frac{c_0}{\Phi_0(3;x)}
=\frac{1/4}{x^3+b_1x^2+b_2x+b_3}$ such that $b_1, b_2$ and $b_3$ satisfy the following condition
\begin{align*}
\max_{b_1,b_2,b_3}\mathrm{R}\left(\mathrm{MC}_0(x)
-\mathrm{MC}_0(x+1)-\frac{12 x^4-1}{(4 x^4+1  )^2}\right).
\end{align*}
By \eqref{Ek-Mathematica-Program}~(also see example 3 below), we obtain that $b_1=-\frac 32, b_2=\frac 54, b_3=-\frac 38$. Then we
take $\mathrm{MC}_1(x)=\frac{\kappa_1}{x+\lambda_1}$ such that $\kappa_1, \lambda_1$ satisfy the following condition
\begin{align*}
\max_{\kappa,\lambda}\mathrm{R}\left(\frac{c_0}{\Phi_0(\nu;x)
+\frac{\kappa}{x+\lambda}}
-
\frac{c_0}{\Phi_0(\nu;x+1)+\frac{\kappa}{x+1+\lambda}}-\frac{12 x^4-1}{(4 x^4+1  )^2}\right).
\end{align*}
By using of \eqref{Ek-Mathematica-Program} again, we find that $\kappa_1=\frac{1}{16}$ and $\lambda_1=-\frac 12$.
Finally, we check directly that the following finite continued fraction
\begin{align}
\frac{\frac 14}{x^3 -\frac 32 x^2 +\frac 54 x -\frac 38
+ \frac{\frac{1}{16}}{x -\frac 12}}
\end{align}
is a solution of the equation~\eqref{Example-1}. For this example, it suffices for us to work for the initial-correction and the first-correction. We can also use \textbf{Mathematica} command ``RSolve" to verify it again. Further, we note that after simplifying the solution above, the final closed form solution $ y(x)=\frac{-1 + 2 x}{2 (1 - 2 x + 2 x^2)^2}$ satisfies \eqref{Example-1} for all real $x$. Lastly, by the telescoping method, we find that
\begin{align}
\sum_{n=1}^{\infty}\frac{12 n^4-1}{(4 n^4+1  )^2}=y(1)=\frac 12.
\end{align}

\bigskip

\noindent{\textbf{Example 2}} 
It seems it impossible to treat this problem by ``RSolve" command of \textbf{Mathematica} software, because its a very huge of computations. By making use of our method, we can find that the following function
\begin{align}
\frac{3/16}{x^7 -\frac 72x^6 + \frac{91}{12}x^5 -\frac{245}{24}x^4 + \frac{1561}{144}x^3 -\frac{749}{
96}x^2 +\frac{9451}{1728}x -\frac{5845}{3456} +\mathrm{MC}_5(x) } \end{align}
is a solution of the equation
\begin{align}
y(x)-y(x+1)- \frac{1 - 480 x^4 + 8736 x^8 - 21504 x^{12} + 5376 x^{16}}{(1 + 4 x^4)^6}=0,
\end{align}
where
\begin{align*}
\mathrm{MC}_5(x)=\K_{j=1}^{5}\left(\frac{\kappa_j}{x-\frac 12}\right),\quad \kappa_1 = \frac{41041}{20736}, \kappa_2 = -\frac{1024}{1353}, \kappa_3 = \frac{243}{41041}, \kappa_4 = -\frac{451}{4368}, \kappa_5 = \frac{1}{48}.
\end{align*}
Certainly, the final simple closed form solution may be simplified as $$\frac{(-1 + 2 x) (-1 - 2 x + 2 x^2) (1 - 6 x + 6 x^2)}{2 (1 - 2 x +2 x^2)^6}.$$
Hence, in a certain sense, our method may be viewed
as a supplement
of Gosper's algorithm~(see Chapter 5 of Petkovsek et al.~\cite{PWZ}, and some Exercises: 1~(c), 3~(b), 3~(f) in~\cite[p.~95\nobreakdash--96]{PWZ}), as our approach does not need to do factors
of polynomials, all things is that we only solve several linear
equations. We shall give some more examples in Sec.~5\nobreakdash--7 below.

\bigskip

\noindent \textbf{(Step 4)} Based on Step 2, construct a formal continued fraction solution of of the difference equation~\eqref{Difference Equation}, then propose a reasonable conjecture. For instance, the continued fraction representation for many mathematical constants such as Catalan constant, $\zeta(2)$, Ap\'ery number $\zeta(3)$, the Mathieu series, the alternating Mathieu series, etc. can be guessed by our method, some of them are new. Lastly, with the help of continued fraction theory and  hypergeometric series, etc, we
try to prove and \emph{extend} the assertion.
\bigskip

\noindent \textbf{(Step 5)} Based on Step 3, one may construct a simple closed form solution of some other equation~\eqref{Difference Equation} by making use of the theory of the linear difference equation, for example, such as the following two properties:

\noindent\textbf{(i)} Let $f_1(x)$ and $f_2(x)$ be two rational function in $x$. If $y_j(x)~(j=1,2)$ is a simple closed form solution of the equation
\begin{align}
y_j(x)-U(x)~y_j(x+1)=f_j(x),
\end{align}
respectively, then for all $c_1,c_2\in \mathbb{R}$, $c_1y_1(x)+c_2y_2(x)$ is also a simple closed form solution of the equation
\begin{align}
y(x)-U(x)~y(x+1)=c_1f_1(x)+c_2f_2(x).
\end{align}

\noindent\textbf{(ii)} Given a rational function $f(x)$. If the equation
\begin{align}
y(x)-c~y(x+1)=f(x),
\end{align}
has a simple closed form solution, so is
\begin{align}
y(x)-c~y(x+1)=\frac{d^k}{dx^k}f(x).
\end{align}

\begin{rem}
In fact, Example 1 and 2 were constructed by property (ii) and the \eqref{Difference Equation-6} in Theorem 9 with $(p,q,s,r)=(1,0,0,1/4)$. All calculations in this work were performed by using of \textbf{Mathematica} version 8.0.
\end{rem}
\bigskip

Now we use the Mathieu series to illustrate how to guess its continued fraction representation.

\noindent{\textbf{Example 3}} Let us consider the equation
\begin{align}
y(x)-y(x+1)-\frac{2x}{(x^2+r^2)^2}=0.\label{Mathieu-DF}
\end{align}
As $\mathrm{R}\left(\frac{2x}{(x^2+r^2)^2}\right)=3$, it is not hard to see that we should choose the initial-correction $\mathrm{MC}_0(x)$ in the form $\mathrm{MC}_0(x)=\frac{c_0}{x^2+d_1 x+d_2}$. We manipulate \textbf{Mathematica} software to expand $E_0(x)$ into a power series in terms of $1/x$
\begin{align}
E_0(x)=&\mathrm{MC}_0(x)-\mathrm{MC}_0(x+1)-\frac{2x}{(x^2+r^2)^2}\\
=&\frac{-2 + 2 c_0}{x^3} - \frac{
 3 (c_0 + c_0 d_1)}{x^4} +\frac{4 c_0 + 6 c_0 d_1 + 4 c_0 d_1^2 - 4 c_0 d_2 + 4 r^2}{x^5}+O\left(\frac{1}{x^6}\right).\nonumber
\end{align}
We enforce the first three coefficients to be zero, and find
\begin{align}
c_0=1,\quad d_1=-1,\quad d_2=\frac{1 + 2 r^2}{2}.
\end{align}
Note that the solution is unique~! By \textbf{Mathematica} software again, one may check that $\mathrm{R}\left(E_0(x)\right)=7$. Repeating the process above several times, one observes that the $k$th-correction $\mathrm{MC}_k(x)$ is a \emph{Type-II} and the $\mathrm{MC}$-point
$\omega=-\frac 12$. As the detail is quite similar to the initial-correction, here we only list final computation results as follows
\begin{align}
\mathrm{MC}_k(x)=\K_{j=1}^{k}\left(\frac{\kappa_j}
{(x+\omega)^2+\lambda_j}\right),
\end{align}
where
\begin{align*}
&\kappa_1=-\frac{1}{12} (1 + 4 r^2),\quad \lambda_1=\frac{5 + 4 r^2}{4},\\
&\kappa_2=-\frac{16}{15} (1 + r^2),\quad \lambda_2=\frac{13 + 4 r^2}{4},\\
&\kappa_3=-\frac{81}{140} (9 + 4 r^2),\quad \lambda_3=\frac{25 + 4 r^2}{4},\\
&\kappa_4=-\frac{256}{63} (4 + r^2),\quad \lambda_4=\frac{41 + 4 r^2}{4},\\
&\kappa_5=-\frac{625}{396} (25 + 4 r^2),\quad \lambda_5=\frac{61 + 4 r^2}{4},\\
&\kappa_6=-\frac{1296}{143} (9 + r^2),\quad \lambda_6=\frac{85 + 4 r^2}{4},\\
&\kappa_7=-\frac{2401}{780} (49 + 4 r^2),\quad \lambda_7=\frac{113 + 4 r^2}{4}.
\end{align*}
Just as did in Sec.~8 of \cite{CTZ}, by careful data analysis and further checking, we may guess that the following formal continued fraction
\begin{align}
\frac{1}{\left(x-\frac 12\right)^2+\frac 14\left(1+4r^2\right)+\K_{n=1}^{\infty}
\left(\frac{\kappa_n}{\left(x-\frac 12\right)^2+\lambda_n}\right)}
\end{align}
should be a solution of the equation \eqref{Mathieu-DF}, where
\begin{align}
\kappa_n= -\frac{n^4\left(n^2 + 4 r^2\right)}{4 (2 n - 1) (2 n + 1)},\quad\lambda_n=\frac 14 (2 n^2 + 2 n + 1 + 4 r^2).
\end{align}

Finally, applying the conjecture above, \eqref{Mathieu-Series} and the telescoping method, we could conjecture further a continued fraction formula for the Mathieu series, which was already proved in~\cite{CTZ1}. Also see Example 8 in Sec.~7 below.
\bigskip

As a briefly summary of this section, we stress that the order of priority for our method is as follows:
The best situation is to look for a \emph{simple closed form} solution, which is a finite continued fraction; The next one is to find a continued fraction solution; The third is to find a linear combination solution
of several continued fractions; The last one is
to find a finite continued fraction approximation solution as you want, some examples, see~\cite{CY}.

On one hand, in order to determine all the related coefficients, we often use an appropriate symbolic computation software, which needs a huge of computations. On the other hand, the exact expression at each occurrence also takes a lot of space. Moreover, in order to guess the continued fraction formula, we have to do a lot of additional works. All theorems in Sec.~5 to 7 are built on experimental results described in this section. Hence, we shall focus on the rigorous proof of all conjectures, and omit the related details for guessing these formulas. Readers interested  may refer to Sec.~6 and 8 in reference~\cite{CTZ}.

\bigskip
\section{Some preliminary lemmas }
In this section, we shall prepare some lemmas for later use. The main lemmas are two of Ramanujan's continued fraction formulas involving the quotient of the gamma functions~(see~\cite{BLW,Ram,Ram1}).

\begin{lem} Let $x,m,$ and $n$ be complex. We define
\begin{align}
Q=Q(x,m,n):=\frac{\Gamma\left(\frac 12(x+m-n+1)\right)\Gamma\left(\frac 12(x-m+n+1)\right)}{\Gamma\left(\frac 12(x+m+n+1)\right)\Gamma\left(\frac 12(x-m-n+1)\right)}.
\end{align}
If either $m$ or $n$ is an integer or if $\Re x>0$, then
\begin{align}
\frac{1-Q}{1+Q}=\frac{mn}{x+\K_{j=1}^{\infty}
\left(\frac{(m^2-j^2)(n^2-j^2)}{
(2j+1)x}\right)}.
\end{align}
\end{lem}
\proof This is Entry 33 of Berndt~\cite[p.~155]{Ber}. For its research history, see p.~156 in~\cite{Ber}.

\begin{lem} Let $x,l,m$, and $n$ denote complex numbers. We define
\begin{align}
P=&P(x,l,m,n)\\
=:&\frac{\Gamma\left(\frac 12(x+l+m+n+1)\right)\Gamma\left(\frac 12(x+l-m-n+1)\right)\Gamma\left(\frac 12(x-l+m-n+1)\right)\Gamma\left(\frac 12(x-l-m+n+1)\right)}{\Gamma\left(\frac 12(x-l-m-n+1)\right)\Gamma\left(\frac 12(x-l+m+n+1)\right)\Gamma\left(\frac 12(x+l-m+n+1)\right)\Gamma\left(\frac 12(x+l+m-n+1)\right)}.\nonumber
\end{align}
Then if either $l,m$, or $n$ is an integer or if $\Re x>0$,
\end{lem}
\begin{align}
\frac{1-P}{1+P}=\frac{2lmn}{x^2-l^2-m^2-n^2+1+\K_{j=1}^{\infty}
\left(
\frac{4(l^2-j^2)(m^2-j^2)(n^2-j^2)}{(2j+1)\left(x^2-l^2-m^2-n^2
+2j^2+2j+1\right)}\right)
}.
\end{align}
\proof This is Entry 35 of B. C. Berndt~\cite[p.~157]{Ber},  which was claimed first by Ramanujan~\cite{Ram1}. The first published proof was provided by Watson~\cite{Wat}. For the full proof of
Lemma 2, we refer the reader to L. Lorentzen's paper~\cite{Jac}.\qed

\begin{lem}
$b_0+\K\left(a_n/b_n\right)\approx d_0+\K\left(c_n/d_n\right)$ if and only if there exists a sequence $\{r_n\}$ of complex numbers with $r_0=1,r_n\neq 0$ for all $n\in \mathbb{N}$, such that
\begin{align}
d_0=b_0,\quad c_n=r_{n-1}r_n a_n,\quad d_n=r_n b_n\quad \mbox{for all
 $n\in \mathbb{N}$.}
\end{align}
\end{lem}
\proof See Theorem 9 in L. Lorentzen, H. Waadeland~\cite[p.~73]{LW} .\qed

\section{The rational series $\sum\frac{1}{n^2+an+b}$}
Let $a,b\in\mathbb{R}$. Given the following infinite series
\begin{align}
S(a,b)=\sum_{n=n_0}^{\infty}\frac{1}{n^2+an+b},
\end{align}
where $n_0$ is a suitable non-negative integer such that
$n^2+an+b\neq 0$
for all integer $n\ge n_0$. It is a natural question when the coefficient
$a$ and $b$ satisfy the conditions, one can get a simple closed
form for the sum of $S(a,b)$. It seems that the results in this section are not very new. However, our purpose is to treat it in a unified setting.
Firstly, we shall study the following difference equation
\begin{align}
y(x)-y(x+1)=\frac{1}{x^2+a x +b}.
\end{align}

\begin{thm} Let $a, b\in \mathbb{R}$, and the formal continued fraction $F(a,b;x)$ or shortly $F(x)$, be defined by
\begin{align}
F(a,b;x):=\frac{1}{x+\omega+\K_{n=1}^{\infty}
\left(\frac{\kappa_n}{x+\omega}\right)
},\label{F-a-b-x-D}
\end{align}
where
\begin{align}
\omega=\frac{a-1}{2},\quad\kappa_n=\frac{n^2(n^2 + 4 b - a^2)}{4 (2 n - 1) (2 n + 1)}.\label{omega-1}
\end{align}
We assume that $x\notin\{q+\alpha: q\in\mathbb{Z}\backslash\mathbb{N},~ {\alpha}^2+a{\alpha} +b=0, ~\alpha\in\mathbb{C} \}$. If either $\sqrt{a^2-4b}\in\mathbb{N}$  or  $\Re x>-\omega$, then
\begin{align}
F(x)-F(x+1)=
\frac{1}{x^2+a x+b}.\label{Difference Equation-1}
\end{align}
\end{thm}

\proof We shall discuss the following three cases.
\bigskip

\noindent \emph{(Case 1)}
Assume $b<\frac{a^2}{4}$.
By Lemma 1 with
$(x,m)=\left(2(x+\omega),\sqrt{a^2-4b}\right)$, under the conditions of Theorem 1 we have
\begin{align}
\frac{1-Q}{1+Q}=\frac{ n\sqrt{a^2-4b}}{2(x+\omega)+\K_{j=1}^{\infty}\left(
\frac{(a^2-4b-j^2)(n^2-j^2)}{2(2j+1)\left(x+\omega\right)}\right)}.
\end{align}
Dividing both sides by $ n\sqrt{a^2-4b}$ and letting $n$ tend to zero, on the right side, we deduce that
\begin{align}
\frac{1}{2(x+\omega)+\K_{j=1}^{\infty}\left(
\frac{j^2(j^2+4b-a^2)}{2(2j+1)\left(x+\omega\right)}\right)}.
\end{align}
The following classical representation is well-known~(e.g., see~\cite[Eq.~6.3.16, p.~259]{AS})
\begin{align}
\psi(z+1)=-\gamma+\sum_{k=1}^{\infty}\left(\frac 1k-\frac{1}{k+z}\right),\quad (z\neq -1,-2,-3,\ldots),\label{PolyGamma-Exprssion}
\end{align}
where $\gamma$ denotes Euler-Mascheroni constant. In the sequel, we shall use this formula several times, usually without comment.

On the other hand, from the definition of $Q$, we see easily that $\lim_{n\rightarrow 0}Q=1$. A direct calculation with the use of L'Hospital's rule gives
\begin{align}
&\frac{1}{\sqrt{a^2-4b}}\lim_{n\rightarrow 0}\frac{1-Q}{n(1+Q)}
=\frac{1}{\sqrt{a^2-4b}}\lim_{n\rightarrow 0}\frac{1}{1+Q}
\lim_{n\rightarrow 0}\frac{1-Q}{n}
=\frac{1}{2\sqrt{a^2-4b}}\lim_{n\rightarrow 0}\frac {\partial}{\partial n}(1-Q)
\\
=&\frac{1}{2\sqrt{a^2-4b}}
\left\{-\psi\left(x+\frac{a-\sqrt{a^2-4b}}{2}\right)
+\psi\left(x+\frac{a+\sqrt{a^2-4b}}{2}\right)\right\}\nonumber\\
=&\frac{1}{2\sqrt{a^2-4b}}\sum_{k=0}^{\infty}\left(\frac{1}
{k+x+\frac{a-\sqrt{a^2-4b}}{2}}-\frac{1}
{k+x+\frac{a+\sqrt{a^2-4b}}{2}}\right)\nonumber\\
=&\frac{1}{2}\sum_{k=0}^{\infty}\frac{1}{(k+x+\frac a2)^2-\frac{a^2-4b}{4}}.\nonumber
\end{align}
Hence,
\begin{align*}
\sum_{k=0}^{\infty}\frac{1}{(k+x+\frac a2)^2-\frac{a^2-4b}{4}}
=\frac{2}{2(x+\omega)+\K_{j=1}^{\infty}\left(
\frac{j^2(j^2+4b-a^2)}{2(2j+1)\left(x+\omega\right)}\right)}
=F(x),\nonumber
\end{align*}
where we used Lemma 3 in the last equality.
Under the conditions of Theorem 1, it is not difficult to prove that
\begin{align}
F(x)-F(x+1)=&\sum_{k=0}^{\infty}\frac{1}{(k+x+\frac a2)^2-\frac{a^2-4b}{4}}-\sum_{k=0}^{\infty}\frac{1}{(k+x+1+\frac a2)^2-\frac{a^2-4b}{4}}\\
=&\frac{1}{(x+\frac a2)^2-\frac{a^2-4b}{4}}
=\frac{1}{x^2+a x+b}.\nonumber
\end{align}
Hence the identity~\eqref{Difference Equation-1} is true in this case.
\bigskip

\noindent \emph{(Case 2)} Suppose $b>\frac{a^2}{4}$,  it follows from Lemma 1 with $(x,m)=\left(2(x+\omega),\sqrt{4b-a^2}~i\right)$ that
\begin{align}
\frac{1-Q}{1+Q}=\frac{ n\sqrt{4b-a^2}~i}{2(x+\omega)+\K_{j=1}^{\infty}\left(
\frac{(a^2-4b-j^2)(n^2-j^2)}{2(2j+1)\left(x+\omega\right)}
\right)}.
\end{align}
By using of L'Hospital's rule, we deduce that
\begin{align}
&\frac{1}{\sqrt{4b-a^2}~i}\lim_{n\rightarrow 0}\frac{1-Q}{n(1+Q)}\\
=&\frac{1}{\sqrt{4b-a^2}~i}\lim_{n\rightarrow 0}\frac{1}{1+Q}
\lim_{n\rightarrow 0}\frac{1-Q}{n}
=\frac{1}{2\sqrt{4b-a^2}~i}\lim_{n\rightarrow 0}\frac {\partial}{\partial n}(1-Q)\nonumber
\\
=&\frac{1}{2\sqrt{4b-a^2}~i}
\left(-\psi\left(x+\frac{a-\sqrt{4b-a^2}~i}{2}\right)
+\psi\left\{x+\frac{a+\sqrt{4b-a^2}~i}{2}\right)\right\}\nonumber\\
=&\frac{1}{2\sqrt{4b-a^2}~i}\sum_{k=0}^{\infty}\left(\frac{1}
{k+x+\frac{a-\sqrt{4b-a^2}~i}{2}}-\frac{1}
{k+x+\frac{a+\sqrt{4b-a^2}~i}{2}}\right)\nonumber\\
=&\frac{1}{2}\sum_{k=0}^{\infty}\frac{1}{(k+x+\frac a2)^2+\frac{4b-a^2}{4}}.\nonumber
\end{align}
Quite similarly to Case 1, Theorem 1 holds true in Case 2.
\bigskip

\noindent \emph{(Case 3)} $b=\frac{a^2}{4}$. In this case, it follows from Lemma 1 with $x=2(x+\omega)$ that
\begin{align}
\lim_{m\rightarrow 0}\frac 1m\lim_{n\rightarrow 0}
\frac{1-Q}{n(1+Q)}=&\lim_{m\rightarrow 0}\lim_{n\rightarrow 0}
\frac{1}{2(x+\omega)+\K_{j=1}^{\infty}\left(
\frac{(m^2-j^2)(n^2-j^2)}{2(2j+1)\left(x+\omega\right)}\right)}\\
=&\frac{1}{2(x+\omega)+\K_{j=1}^{\infty}\left(
\frac{j^4}{2(2j+1)\left(x+\omega\right)}\right)}=\frac 12 F(x).\nonumber
\end{align}
By making use of L'Hospital's rule twice, we find that
\begin{align}
&\lim_{m\rightarrow 0}\frac 1m\lim_{n\rightarrow 0}
\frac{1-Q}{n(1+Q)}=\lim_{m\rightarrow 0}\frac 1m\left(
\lim_{n\rightarrow 0}\frac{1}{1+Q}\lim_{n\rightarrow 0}
\frac{1-Q}{n}
\right)\\
=&\frac 12\lim_{m\rightarrow 0}\frac 1m\left\{-\psi\left(x+\frac{a-m}{2}\right)
+\psi\left(x+\frac{a+m}{2}\right)\right\}\nonumber\\
=&\frac 12\psi'(x+\frac a2)=\frac 12\sum_{k=1}^{\infty}\frac{1}{(k-1+x+\frac{a}{2})^2}.\nonumber
\end{align}
Applying the similar argument as the proof of Case 1, one has
\begin{align}
F(x)-F(x+1)=&\sum_{k=1}^{\infty}\frac{1}{(k-1+x+\frac{a}{2})^2}
-\sum_{k=1}^{\infty}\frac{1}{(k+x+\frac{a}{2})^2}\label{Theorem-1-Case 3}\\
=&\frac{1}{(x+\frac a2)^2}=\frac{1}{x^2+a x+b}.\nonumber
\end{align}
We remark that combining  Case 1, Case 2 and a limiting process for ~\eqref{Difference Equation-1}~(i.e. let $b$ tend to $\frac{a^2}{4}$), the \eqref{Theorem-1-Case 3} may be proved easily. This completes the proof of Theorem 1.\qed

\begin{thm} With the notations of Theorem 1, let $n_0$ is a non-negative integer such that  $n_0>\left\{\alpha:~\alpha\in\{\frac {-a\pm\sqrt{a^2-4b}}{2}\}~\text{and}~\alpha\in\mathbb{Z}\right\}$. If either $\sqrt{a^2-4b}\in \mathbb{N}$ or $n_0>\frac{-a+1}{2}$, then
\begin{align}
\sum_{n=n_0}^{\infty}\frac{1}{n^2+an+b}=F(n_0).
\end{align}
In particular, if
\begin{align}
b\in \left\{\frac{a^2-k^2}{4}: a\in \mathbb{R}, k\in \mathbb{N}\right\},
\end{align}
and $n_0>\max\left\{\alpha: \alpha\in\{-1,\frac {-a\pm k}{2}\}~\text{and}~\alpha\in\mathbb{Z}\right\}$, then
\begin{align}
\sum_{n=n_0}^{\infty}\frac{1}{n^2+an+b}=
\frac{1}{n_0+\frac{a-1}{2}+\K_{n=1}^{k-1}
\left(\frac{\frac{n^2(n^2 -k^2)}{4 (2 n - 1) (2 n + 1)}}{n_0+\frac{a-1}{2}}\right)
}.\label{Theorem-2-2}
\end{align}
\end{thm}
\proof It follows readily from Theorem 1 and the telescoping method.\qed


\bigskip

\noindent\textbf{Example 4.} We let $k=3$, $b=\frac{a^2-9}{4}$, $\omega=\frac{a-1}{2}$, and
\begin{align}
F_a(x)=\frac{1}{x+\omega+\K_{n=1}^{2}
\left(\frac{\kappa_n}{x+\omega}\right)
}
=\frac{2 (-1 - 12 x + 12 x^2 - 6 a + 12 x a +
3 a^2)}{3 (-1 + 2 x + a) (-3 - 4 x + 4 x^2 -
2 a + 4 x a + a^2)}.\nonumber
\end{align}
We can check directly that if $x\neq \frac{1-a}{2}+q, \frac{-a\pm 3}{2}+q$, $q\in \{0,-1\}$, then
\begin{align}
F_a(x)-F_a(x+1)=\frac{1}{x^2+a x+\frac{a^2-9}{4}}.\label{Example 1-1}
\end{align}

Let $p$ be prime and set $a=\sqrt{p}$, then the following series is an irrational number
\begin{align}
\sum_{n=1}^{\infty}\frac{1}{n^2+\sqrt {p} n+\frac {p-9}{4}}=F_{\sqrt {p}}(1)=\frac{6 (p - 1) - 2 \sqrt{p} (3 p - 19)}{3 (-9 + 10 p - p^2)}.\label{Example 1-2}
\end{align}

For instance, $
F_{\sqrt {2}}(1)=\frac{2(3 + 13 \sqrt{2})}{21}$ is irrational. As a by-product of Theorem 2, one may employ the infinite series in \eqref{Theorem-2-2} to construct many irrational numbers . Moreover, we can check \eqref{Example 1-1} and \eqref{Example 1-2} by applying of \textbf{Mathematica} commands ``RSolve" and ``Sum", respectively.
\bigskip

\noindent\textbf{Example 5.} Let $p,q\in\mathbb{N}$ with $p>q$ and $(p,q)=1$, and $r>0$, then the following rational series
\begin{align}
\sum_{n=1}^{\infty}\frac{1}{(p n+q)^2},\quad \sum_{n=1}^{\infty}\frac{1}{n^2+r^2}
\end{align}
can be written as a continued fraction expansion. For complex $r$, Ramanujan even deduced a exact expression for the second series above, see Entry 24(i) and (ii) in~\cite[Chap.~14, p.~291\nobreakdash--292]{Ber}. From Whittaker and Watson's text~\cite[p.~136]{WW}, one has
\begin{align}
\frac{1}{e^x-1}=\frac 1x-\frac 12+\sum_{m=1}^{\infty}\frac{2x}{x^2+4\pi ^2m^2}.
\end{align}
However, our approach is different from their methods.

\section{The rational series $\sum\frac{1}{Q_3(n)}$}

Let $a,b,c\in\mathbb{R}$ and $Q_3(x)=x^3+a x^2+b x+c$. Similarly to the previous section, we first study the following difference equation
of order one
\begin{align}
y(x)-y(x+1)=\frac{1}{Q_3(x)}.\label{Difference Equation-2}
\end{align}
From the fundamental theorem of algebra, a polynomial of degree three with real coefficients may be expressed as
\begin{align}
Q_3(x)=(x+t)(x^2+r x+s),
\end{align}
here $r,s,t\in\mathbb{R}$. If the discriminant $\Delta=r^2-4s\ge 0$ for the last polynomial of degree two, then we further write it in the form
\begin{align}
Q_3(x)=(x+t)(x+\alpha)(x+\beta),\quad \alpha,\beta\in \mathbb{R}.
\end{align}
If $\alpha=\beta=t$, then it reduces to Entry 32~(iii) in Berndt~\cite{Ber}, also see Case 3 in the proof of Theorem 3 below. Otherwise, without loss of the generality, we assume $\alpha\neq \beta$. In which case, it follows from Theorem 1 that for
$\Re x>\max\{-\frac{\alpha+t-1}{2},-\frac{\beta+t-1}{2}\}$
\begin{align}
&\frac{1}{Q_3(x)}=\frac{1}{(x+t)(x+\alpha)(x+\beta)}
=\frac{1}{\beta-\alpha}\left(
\frac{1}{(x+t)(x+\alpha)}-\frac{1}{(x+t)(x+\beta)}\right)
\label{trivial assertion-1}\\
=&\frac{1}{\beta-\alpha}\left\{(F(t+\alpha,t\alpha;x)-F(t+\alpha,t \alpha;x+1))
-(F(t+\beta,t \beta;x)-F(t+\beta,t \beta;x+1))\right\},\nonumber
\end{align}
where $F(a,b;x)$ is given as~\eqref{F-a-b-x-D}. By Theorem 1, it is not hard to get the following assertion.

\noindent\textbf{Example 6.} Let $\upsilon\in \mathbb{R}$, $k\in \mathbb{Z}\backslash\{0,1\}$ and $d\in \mathbb{Z}\backslash\{0\}$. If $x>\max\{-\upsilon, -\upsilon-d, -\upsilon-k d\}$, then the following equation
\begin{align}
y(x)-y(x+1)=\frac{1}{(x+\upsilon)(x+\upsilon+d)(x+\upsilon+k d)},\label{Example -4}
\end{align}
has a simple closed form solution.

If $\Delta=r^2-4s<0$, to the best knowledge of authors, up to now
very little has been established except in the form of complex function.

A lot of experiment results show that the structure of the continued fraction solution for the equation~\eqref{Difference Equation-2} may  be the \emph{type I} or \emph{type II} or other type according to the various conditions of the parameters $a,b$ and $c$. In this section, we shall apply our method to find new results for two special classes, and then give further remarks.

\subsection{For the case of $c=\frac{-2 a^3 + 9 a b}{27}$}

\begin{thm} Let $a, b\in\mathbb{R}$ and $c=\frac{-2 a^3 + 9 a b}{27}$. Let the formal continued fraction $G_1(x)$ be defined by
\begin{align}
G_1(x):=\frac{1/2}{(x+\omega)^2+\frac{3 - 2 a^2 + 6 b}{12}+\K_{n=1}^{\infty}
\left(\frac{\kappa_n}{(x+\omega)^2+\lambda_n}\right)},
\end{align}
where $\omega=\frac{2 a-3}{6}$ and
\begin{align}
\kappa_n=-\frac{n^2\left(-3 n^2 + a^2 - 3 b\right)^2}{6^2 (2 n -1)
(2 n + 1)},\quad
\lambda_n=\frac{3 - 2 a^2 + 6 b}{12} + \frac{ n +  n^2}{2}.
\end{align}
Assume that $x\notin\{q+\alpha: q\in\mathbb{Z}\backslash\mathbb{N},~ {\alpha}^3+a{\alpha}^2 +b\alpha +c=0,~\alpha\in\mathbb{C} \}$. If either $\sqrt{(a^2-3 b)/3})\in\mathbb{N} $ or $\Re x>-\omega$, then
\begin{align}
G_1(x)-G_1(x+1)=\frac{1}{x^3 + a x^2 + b x + c}.\label{Difference Equation-3}
\end{align}
\end{thm}
\proof It is not difficult to verify that
\begin{align}
x^3 + a x^2 + b x + c
=\frac{1}{27} (a + 3 x) (-2 a^2 + 9 b + 6 a x + 9 x^2),
\end{align}
and the last polynomial of degree $2$ above has the discriminant
\begin{align}
\Delta=(6a)^2-4\cdot 9(-2 a^2 + 9 b)=108 (a^2 - 3 b)\begin{cases}
\geq 0,\quad \mbox{if $b\leq \frac{a^2}{3}$},\\
<0,\quad \mbox{otherwise.}
\end{cases}
\end{align}
We shall consider three cases.

\noindent \emph{(Case 1)} $a^2-3 b>0$. Applying Lemma 2 with $(x,l,m)=(2(x+\omega),\sqrt{(a^2-3 b)/3},\sqrt{(a^2-3 b)/3})$ and dividing both sides by $2 n(a^2-3 b)/3$, under the conditions of Theorem 5 we obtain that

\begin{align}
\frac{3}{2(a^2-3 b)}\frac{1-P}{n(1+P)}
=\frac{1}{4(x+\omega)^2-n^2-\frac{2(a^2-3 b)}{3}+1+\K_{j=1}^{\infty}\left(
\frac{4\left(\frac{a^2-3b}{3}-j^2\right)^2(n^2-j^2)}{(2j+1)
\left(4(x+\omega)^2-n^2-\frac{2(a^2-3 b)}{3}
+2j^2+2j+1\right)}\right)}.
\end{align}
Now let $n$ tend to zero. On the right side, we arrive at
\begin{align}
\frac{1}{4(x+\omega)^2-\frac{2(a^2-3 b)}{3}+1+\K_{j=1}^{\infty}\left(
\frac{-4j^2\left(\frac{a^2-3b}{3}-j^2\right)^2}{(2j+1)
\left(4(x+\omega)^2-\frac{2(a^2-3 b)}{3}
+2j^2+2j+1\right)}\right)}
=\frac 12 G_1(x).
\end{align}
On the other hand, from the definition of $P$, we observe easily that $\lim_{n\rightarrow 0}P=1$. A direct calculation with the use of L'Hospital's rule gives
\begin{align}
&\lim_{n\rightarrow 0}\frac{1-P}{n(1+P)}=\lim_{n\rightarrow 0}
\frac{1}{1+P}\lim_{n\rightarrow 0}\frac{1-P}{n}
\lim_{n\rightarrow 0}=\frac 12\lim_{n\rightarrow 0}
\frac {\partial}{\partial n}(1-P)
\\
=&\psi\left( x+\frac a3\right)-\frac 12\psi\left(x+\frac a3+\sqrt{\frac{a^2-3b}{3}}\right)-\frac 12\psi\left(x+\frac a3-\sqrt{\frac{a^2-3b}{3}}\right)\nonumber\\
=&\frac 12\sum_{k=0}^{\infty}\left(-\frac{2}{k+x+a/3}+\frac{1}{k+x+a/3+
\sqrt{\frac{a^2-3b}{3}}}+\frac{1}{k+x+a/3-
\sqrt{\frac{a^2-3b}{3}}}
\right)\nonumber\\
=&\frac{a^2-3b}{3}\sum_{k=0}^{\infty}\frac{1}{(k + x+\frac a3)\left((k + x+\frac a3)^2+3b-a^2\right)}.\nonumber
\end{align}
Therefore,
\begin{align}
&G_1(x)-G_1(x+1)\\
=&\sum_{k=0}^{\infty}\frac{1}{( k +  x+\frac a3)\left(( k +  x+
\frac a3)^2-a^2+3b\right)}-\sum_{k=0}^{\infty}\frac{1}{(k +  x+1+\frac a3)\left(( k +  x+1+\frac a3)^2-a^2+3b\right)}\nonumber\\
=&\frac{1}{x^3 + a x^2 + b x + c}.\nonumber
\end{align}
This completes the proof of Theorem 3 in this case.
\bigskip

\noindent \emph{(Case 2)} $a^2-3 b<0$. Applying Lemma 2 with $(x,l,m)=(2(x+\omega),\sqrt{(3 b-a^2)/3}~i,\sqrt{(3 b-a^2)/3}~i)$ and dividing both sides by $2 n(a^2-3 b)/3$, similarly to Case 1, we also have

\begin{align}
&\lim_{n\rightarrow 0}\frac{1-P}{n(1+P)}=\lim_{n\rightarrow 0}
\frac{1}{1+P}\lim_{n\rightarrow 0}\frac{1-P}{n}
\lim_{n\rightarrow 0}=\frac 12\lim_{n\rightarrow 0}
\frac {\partial}{\partial n}(1-P)\\
=&\psi\left( x+\frac a3\right)-\frac 12\psi\left(x+\frac a3+\sqrt{\frac{3b-a^2}{3}}~i\right)-\frac 12\psi\left(x+\frac a3-\sqrt{\frac{3b-a^2}{3}}~i\right)\nonumber\\
=&\frac 12\sum_{k=0}^{\infty}\left(-\frac{2}{k+x+a/3}+\frac{1}{k+x+a/3+
\sqrt{\frac{3b-a^2}{3}}~i}+\frac{1}{k+x+a/3-
\sqrt{\frac{3b-a^2}{3}}~i}
\right)\nonumber\\
=&\frac{a^2-3b}{3}\sum_{k=0}^{\infty}\frac{1}{(k + x+\frac a3)\left((k + x+\frac a3)^2+3b-a^2\right)}.\nonumber
\end{align}
Hence, Theorem 3 is true in Case 2.

\bigskip

\noindent \emph{(Case 3)} $a^2-3 b=0$. Applying Lemma 2 with $(x,l,m)=(2(x+\omega),n, n)$,  dividing both sides by $2n^3$,
and employing L'Hospital's rule three times, then
\begin{align}
\lim_{n\rightarrow 0}\frac{1-P}{n^3(1+P)}=&
\lim_{n\rightarrow 0}
\frac{1}{1+P}\lim_{n\rightarrow 0}\frac{1-P}{n^3}
\lim_{n\rightarrow 0}=\frac 12\lim_{n\rightarrow 0}\frac{1}{3n^2}
\frac {\partial}{\partial n}(1-P)
\\
=&-\frac{1}{2}\psi''\left(x+
\frac a3\right)=\sum_{k=0}^{\infty}\frac{1}{\left(k+x+
\frac a3\right)^3}.\nonumber
\end{align}
Now Theorem 3 follows from the trivial equality
\begin{align*}
x^3 + a x^2 + b x + c=\left(x+
\frac a3\right)^3.
\end{align*}
Lastly, combining three cases above will finish the proof of Theorem 3.\qed
\begin{thm} With the notations of Theorem 3, let $n_1$ be a non-negative integer such that $n_1>\max\left\{\alpha:~\alpha\in\{-\frac a3, -\frac a3\pm\sqrt{\frac{a^2-3b}{3}}\}~\text{and}~\alpha\in \mathbb{Z}\right\}$. If
either $\sqrt{(a^2-3 b)/3}\in\mathbb{N} $ or $n_1>-\omega$ , then
\begin{align}
\sum_{n=n_1}^{\infty}\frac{1}{n^3+an^2+bn +\frac{-2 a^3 + 9 a b}{27}}=G_1(n_1).
\end{align}
In particular, if
\begin{align}
b\in\mathfrak{D}_1=\{\frac{a^2}{3}-k^2: k\in \mathbb{N}, a\in\mathbb{R}\},
\end{align}
and $n_1>\max\left\{\alpha:~\alpha\in\{-1, -\frac a3, -\frac a3\pm k\}~\text{and}~\alpha\in \mathbb{Z}\right\}$, then
\begin{align}
\sum_{n=n_1}^{\infty}\frac{1}{n^3+an^2+bn +\frac{-2 a^3 + 9 a b}{27}}=\frac{1/2}{(n_1+\omega)^2+\frac{1-2k^2}{4}+\K_{n=1}^{k-1}
\left(\frac{-\frac{n^2(k^2-n^2)^2}{4(2n-1)(2n+1)}}
{(n_1+\omega)^2+\frac{1-2k^2}{4}+\frac{n+n^2}{2}}\right)}.
\end{align}
\end{thm}
\proof It follows from Theorem 3 and the telescoping method.\qed

\subsection{For the infinite series
$\sum\frac{1}{(n+u)\left(n+\frac{u+v}{2}\right)
(n+v)}$}

\begin{thm} Let $u, v\in\mathbb{R}$. Let the formal continued fraction $G_2(x)$ be defined by
\begin{align}
G_2(x):=\frac{1/2}{(x+\omega)^2+\lambda_0+\K_{n=1}^{\infty}
\left(\frac{\kappa_n}{(x+\omega)^2+\lambda_n}\right)},
\end{align}
where
\begin{align}
&\omega=-\frac 12+u-v,\quad\lambda_n=\frac{2n^2 + 2n + 1}{4}-\frac{(u-v)^2}{8},\\
&\kappa_n=-\frac{n^2 \left(- n^2 + \left(\frac{ u-v}{2}\right)^2\right)^2 }{4 (2 n - 1) (2 n + 1)}.
\end{align}
Let $x\notin\left\{\alpha +q: q\in\mathbb{Z}\backslash\mathbb{N}, \alpha\in\{-u,-v, -\frac{u+v}{2} \}\right\}$. If either $\frac{ u-v}{2}\in\mathbb{Z}\backslash\{0\}$ or $\Re x>-\omega$,   then
\begin{align}
G_2(x)-G_2(x+1)=
\frac{1}{(x+u)\left(x+\frac{u+v}{2}\right)
\left(x+v\right)} .\label{Difference Equation-4}
\end{align}
\end{thm}
\proof In what follows, we always assume $u\neq v$, otherwise it turns to Case 3 in the proof of Theorem 3. We set $t=\frac{(u - v)^2}{8}$. Applying Lemma 2 with $(x,l,m)=(2(x+\omega),\frac{u-v}{2},\frac{u-v}{2})$ and dividing both sides by $\frac{n(u-v)^2}{2}$, under the conditions of Theorem 5 we get that
\begin{align}
\frac{2}{(u-v)^2}\frac{1-P}{n(1+P)}
=\frac{1}{(2x+2\omega)^2-n^2-4t+1+\K_{j=1}^{\infty}\left(
\frac{4\left(\left(\frac{u-v}{2}\right)^2-j^2\right)^2(n^2-j^2)}{(2j+1)
\left((2x+2\omega)^2-n^2-4t
+2j^2+2j+1\right)}\right)}.\nonumber
\end{align}
Now let $n$ tend to zero. On the right side, we arrive at
\begin{align}
\frac{1}{(2x+2\omega)^2-4t+1+\K_{j=1}^{\infty}\left(
\frac{-4j^2\left(\left(\frac{u-v}{2}\right)^2-j^2\right)^2}{(2j+1)
\left((2x+2\omega)^2-4t
+2j^2+2j+1\right)}\right)}=\frac 12 G_2(x).
\end{align}
On the other hand,
a direct calculation with the use of L'Hospital's rule deduces
\begin{align}
&\lim_{n\rightarrow 0}\frac{1-P}{n(1+P)}=\lim_{n\rightarrow 0}
\frac{1}{1+P}\lim_{n\rightarrow 0}\frac{1-P}{n}
\lim_{n\rightarrow 0}=\frac 12\lim_{n\rightarrow 0}
\frac {\partial}{\partial n}(1-P)\\
=&-\frac 12\psi\left( x+u\right)-\frac 12\psi\left(x+v\right)+\psi\left(x+
\frac{u+v}{2}\right)\nonumber\\
=&\frac 12\sum_{k=0}^{\infty}\left(\frac{1}{k+x+u}+\frac{1}{k+x+v}
-\frac{2}{k+x+\frac{u+v}{2}}
\right)\nonumber\\
=&\frac{(u-v)^2}{2}\sum_{k=0}^{\infty}\frac{1}{(k+x+u)
\left(k+x+\frac{u+v}{2}\right)\left(k+x+v\right)
}.\nonumber
\end{align}
Thus,
\begin{align}
&G_2(x)-G_2(x+1)\\
=&\sum_{k=0}^{\infty}\frac{1}{(k+x+u)
\left(k+x+\frac{u+v}{2}\right)\left(k+x+v\right)}
\nonumber\\
&-\sum_{k=0}^{\infty}\frac{1}{(k+1+x+u)
\left(k+1+x+\frac{u+v}{2}\right)\left(k+1+x+v\right)}
\nonumber\\
=&\frac{1}{(x+u)
\left(x+\frac{u+v}{2}\right)\left(x+v\right)}.\nonumber
\end{align}
This completes the proof of Theorem 5.\qed

\begin{thm} With the notations of Theorem 5, let $n_2$ be a non-negative integer such that  $n_2>\max\left\{\alpha:~\alpha\in\{-u, -v,-(u+v)/2\}~\text{and}~\alpha\in\mathbb{Z}\right\}$. If either $\frac{u-v}{2}\in\mathbb{Z}\backslash\{0\}$ or $n_2>\frac 12-u+v$, then
\begin{align}
\sum_{n=n_2}^{\infty}\frac{1}{(n+u)\left(n+\frac{u+v}{2}\right)
\left(n+v\right)}=G_2(n_2).
\end{align}
In particular, if
\begin{align}
u\in\mathfrak{D}_2=\{v\pm 2k: k\in \mathbb{N}, w\in\mathbb{R}\},
\end{align}
and $n_2>\max\{\alpha:~\alpha\in\{-1,-v, -v\mp 2k\}~\text{and}~\alpha\in\mathbb{Z}\}$, then
\begin{align}
&\sum_{n=n_2}^{\infty}\frac{1}{(n+u)\left(n+\frac{u+v}{2}\right)
\left(n+v\right)}=\sum_{n=n_2}^{\infty}\frac{1}{(n+v\pm 2k)\left(n+v\pm k\right)
\left(n+v\right)}\\
=&\frac{1/2}{(n_2-\frac 12\pm 2k)^2+\frac{1-2k^2}{4}+\K_{n=1}^{k-1}
\left(\frac{-\frac{n^2(k^2-n^2)^2}{4(2n-1)(2n+1)}}{(n_2-\frac 12\pm 2k)^2+\frac{2n^2+2n+1}{4}-\frac{k^2}{2}}\right)}.\nonumber
\end{align}
\end{thm}
\proof It follows from Theorem 5 and the telescoping method at once.\qed

\subsection{Some remarks}

\noindent {\textbf{(1)}}  Except the two cases that the polynomial $x^3+ax^2+bx+c$ satisfies the condition of Theorem 3 or 5, for other cases,
the structure of the continued fraction solution of the equation~\eqref{Difference Equation-2} is not a \emph{Type-II}.
\bigskip

\noindent {\textbf{(2)}} Let $r>0$. Taking $a=0$ and $b=r^2$ in Theorem 4, one may obtain a continued fraction representation for the infinite series
\begin{align}
\sum_{n=1}^{\infty}\frac{1}{n(n^2+r^2)}.
\end{align}
When $r\in\mathbb{Q}\backslash\{0\}$, to the best knowledge of authors, the irrationality of the series above remains unproven.

\bigskip

\noindent {\textbf{(3)}}
If $u,v,w\in \mathbb{R}$ with $u\neq v$ and $w\neq 0$~(In fact, if $u=v$, it may be treated by Theorem 3), one has
\begin{align}
\frac{1}{(x+u)\left((x+v)^2+w^2\right)}=&
\frac{1}{(x+u)(x+v-w~i)(x+v+w~i)}\\
=&\frac {1}{2w~i}\left(
\frac{1}{(x+u)(x+v-w~i)}-\frac{1}{(x+u)(x+v+w~i)}
\right).\nonumber
\end{align}
For the following difference equation
\begin{align}
y(x)-y(x+1)=\frac{1}{(x+u)\left((x+v)^2+w^2\right)},
\label{Complex-Case}
\end{align}
similarly to \eqref{trivial assertion-1},  by making use of \textbf{Mathematica} software, authors have checked that if
$\Re \left(x+\frac{u+v-1}{2}\right)>0$, then
\begin{align}
\frac {1}{2w~i}\left\{F(u+v-w~i,u(v-w~i);x)-F(u+v+w~i,u(v+w~i);x)\right\}
\end{align}
is also a solution of the equation~\eqref{Complex-Case}. Here $F(a,b;x)$ is defined as~\eqref{F-a-b-x-D}. For the summation of rational series by means of polygamma functions, please see Section 6.8 in~\cite[p.~264\nobreakdash--265]{AS}.
\bigskip


\section{Two extended Mathieu series}
For a rational series with the general term $u(n)$ in the form of
$u(n)=P_1(n)/P_4(n)$, where $P_j(x)$ is a polynomial of degree $j$ in $x$ with real coefficients, the question become very complexity. Hence in this section, we shall study only two kind of extended Mathieu series.
\subsection{The rational series $\sum\frac{2n+a}{(n^2+a n+b_1)(n^2+a n+b_2)}$}
In this subsection we shall study first the following difference equation of order one
\begin{align*}
y(x)-y(x+1)=
\frac{2x+a}{(x^2+a x+b_1)(x^2+a x+b_2)},
\end{align*}
and the results may is stated as follows.
\begin{thm}
Let $a,b_1,b_2\in\mathbb{R}$, and the formal continued fraction $H_1(a,b_1,b_2;x)$~(or shortly $H_1(x)$) be defined by
\begin{align}
H_1(a,b_1,b_2;x)=\frac{1}{\left(x+\omega\right)^2
+\frac{1+2b_1+2b_2-a^2}
{4}+\K_{n=1}^{\infty}\left(
\frac{\kappa_n}{\left(x+\omega\right)^2+\lambda_n}\right)
},\label{H1-Def}
\end{align}
where
\begin{align}
&\omega=\frac{a-1}{2}, \quad \kappa_n=\frac{n^2 \left(- (b1 - b2)^2+\left(a^2 - 2 (b1 + b2)\right) n^2 - n^4 \right)}{4 (2 n - 1) (2 n + 1)} ,\\
&\lambda_n=\frac{2 n^2 + 2 n + 1+2b_1+2b_2-a^2}{4}.
\end{align}
We assume $x\notin\{q+\alpha: q\in\mathbb{Z}\backslash\mathbb{N},~ ({\alpha}^2+a{\alpha} +b_1)({\alpha}^2+a{\alpha} +b_2)=0, ~\alpha\in\mathbb{C} \}$. If either one of $\sqrt{\frac{a^2 - 2 (b1 + b2)\pm\sqrt{(a^2-4b_1)(a^2-4b_2)}}{2}}$ is a positive integer, or $\Re x>-\omega$, then
\begin{align}
H_1(x)-H_1(x+1)=
\frac{2x+a}{(x^2+a x+b_1)(x^2+a x+b_2)}.
\end{align}
\end{thm}

\proof In this and next subsection, when $t<0$, we shall use the convention $\sqrt{t}=\sqrt{-t}~i$. It is not difficult to prove that
\begin{align}
(b1 - b2)^2-\left(a^2 - 2 (b1 + b2)\right) n^2 + n^4
=\left(\frac{\beta+\sqrt{\Delta_1}}{2}-n^2\right)
\left(\frac{\beta-\sqrt{\Delta_1}}{2}-n^2\right).
\end{align}
where $\beta=a^2 - 2 (b1 + b2)$ and $\Delta_1=(a^2-4b_1)(a^2-4b_2)$.
Since the proof of Theorem 7 is quite similar to that of Theorem 5 or Theorem 9 below, we only give its outline. We shall discuss the following nine cases.


\noindent \emph{(Case 1)}  $\Delta_1>0$ and $\beta-\sqrt{\Delta_1}>0$.
We take $(x,l,m)=\left(2(x+\omega),\sqrt{\frac{\beta-\sqrt{\Delta_1}}{2}},
\sqrt{\frac{\beta+\sqrt{\Delta_1}}{2}}\right)$ in Lemma 2, then let $n$ tend to zero.

\bigskip

\noindent \emph{(Case 2)}  $\Delta_1>0$,  $\beta+\sqrt{\Delta_1}>0$, and $\beta-\sqrt{\Delta_1}<0$.
We take $(x,l,m)=\left(2(x+\omega),\sqrt{\frac{-\beta+\sqrt{\Delta_1}}{2}}~i,
\sqrt{\frac{\beta+\sqrt{\Delta_1}}{2}}\right)$ in Lemma 2, then let $n$ tend to zero.

\bigskip

\noindent \emph{(Case 3)}  $\Delta_1>0$ and $\beta+\sqrt{\Delta_1}<0$.
We take $(x,l,m)=\left(2(x+\omega),\sqrt{\frac{-\beta+\sqrt{\Delta_1}}{2}}~i,
\sqrt{\frac{-\beta-\sqrt{\Delta_1}}{2}}~i\right)$ in Lemma 2, then let $n$ tend to zero.
\bigskip

\noindent \emph{(Case 4)}  $\Delta_1>0$ and $\beta-\sqrt{\Delta_1}=0$.
We take $(x,l)=\left(2(x+\omega),\sqrt{\frac{\beta+\sqrt{\Delta_1}}{2}}
\right)$ in Lemma 2, then let $m$ and $n$ tend to zero, successively.
\bigskip

\noindent \emph{(Case 5)}  $\Delta_1>0$  and $\beta+\sqrt{\Delta_1}=0$.
We take $(x,l)=\left(2(x+\omega),\sqrt{\frac{-\beta
+\sqrt{\Delta_1}}{2}}~i\right)$ in Lemma 2, then let $m$ and $n$ tend to zero, successively.
\bigskip

\noindent \emph{(Case 6)}  $\Delta_1<0$. We take $(x,l,m)=\left(2(x+\omega),\sqrt{\frac{\beta-\sqrt{-\Delta_1}~i}{2}},
\sqrt{\frac{\beta+\sqrt{-\Delta_1}~i}{2}}\right)$ in Lemma 2, then let $n$ tend to zero.
\bigskip

\noindent \emph{(Case 7)}  $\Delta_1=0$ and $\beta>0$. We take $(x,l,m)=\left(2(x+\omega),\sqrt{\frac{\beta}{2}},
\sqrt{\frac{\beta}{2}}\right)$ in Lemma 2, then let $n$ tend to zero.
\bigskip

\noindent \emph{(Case 8)}  $\Delta_1=0$ and $\beta<0$. We take $(x,l,m)=\left(2(x+\omega),\sqrt{-\frac{\beta}{2}}~i,
\sqrt{-\frac{\beta}{2}}~i\right)$ in Lemma 2, then let $n$ tend to zero.

\bigskip

\noindent \emph{(Case 9)}  $\Delta_1=0$ and $\beta=0$.
In this case, we have $b_1=b_2$. It is same as Case 3 in Theorem 3.
Finally, combining the nine cases above will finish the proof of Theorem 7.\qed

\begin{thm} With the notations of Theorem 7, let $n_1$ be a non-negative integer such that $n_1>\max\limits_{\alpha}\{\alpha:~(\alpha^2+a\alpha+b_1)
(\alpha^2+a\alpha+b_2)=0,\alpha\in\mathbb{Z}\}$. If either
one of $\sqrt{\frac{a^2 - 2 (b1 + b2)\pm\sqrt{(a^2-4b_1)(a^2-4b_2)}}{2}}$ is a positive integer, or $n_1>\frac{-a+1}{2}$, then
\begin{align}
\sum_{n=n_1}^{\infty}
\frac{2n+a}{(n^2+a n+b_1)(n^2+a n+b_2)}=H_1(n_1)\label{AN extended MS-1}.
\end{align}
In particular, if
\begin{align}
a\in\left\{\pm\sqrt{k^2 +2 (b1 + b2)+\frac{(b_1-b_2)^2}{k^2}}\in\mathbb{R}:
~b_1,b_2\in\mathbb{R}, k\in\mathbb{N}
\right\},
\end{align}
and $n_1>\max\limits_{\alpha}\{-1, \alpha:~(\alpha^2+a\alpha+b_1)
(\alpha^2+a\alpha+b_2)=0, \alpha\in\mathbb{Z}\}$, then
\begin{align}
\sum_{n=n_1}^{\infty}
\frac{2n+a}{(n^2+a n+b_1)(n^2+a n+b_2)}=\frac{1}{\left(n_1+\omega\right)^2
+\frac{1+2b_1+2b_2-a^2}
{4}+\K_{n=1}^{k-1}\left(
\frac{\kappa_n}{\left(n_1+\omega\right)^2+\lambda_n}\right)
}\label{AN extended MS-1-1}.
\end{align}
Further, if $b_1=b_2=b$, $b\in \{\frac{a^2-k^2}{4}:~a\in \mathbb{R}, k\in\mathbb{N}\}$ and $n_1>\max\left\{\alpha:~\alpha\in\{-1,\frac{-a\pm k}{2}\}~\text{and}~\alpha\in\mathbb{Z}\right\}$, then
\begin{align}
\sum_{n=n_1}^{\infty}\frac{2n+a}{(n^2+a n+b)^2}=\frac{1}{\left(n_1+\frac{a-1}{2}\right)^2+\frac{1-k^2}
{4}+\K_{n=1}^{k-1}\left(
\frac{\frac{n^4\left(k^2-n^2\right)}{4 (2 n - 1) (2 n + 1)}}{\left(n_1+\frac{a-1}{2}\right)^2+\frac{2 n^2 + 2 n + 1-k^2}{4}}\right)
}.\label{AN extended MS-1-2}
\end{align}
\end{thm}
\proof Applying Theorem 7 and the telescoping method will finish the proof of Theorem 8.\qed
\bigskip

\noindent\textbf{Example 7} Taking $k=2, a=\sqrt{41}/2, b_1=2$ and $b_2=1$ in~\eqref{AN extended MS-1-1}, we find that the function
\begin{align}
\frac{12 - \sqrt{41} - 4 x + 2 \sqrt{41} x + 4 x^2}{35 - 5 \sqrt{41} - 61 x + 12 \sqrt{41} x + 65 x^2 - 6 \sqrt{41} x^2 - 8 x^3 +
4 \sqrt{41} x^3 + 4 x^4}
\end{align}
is a simple closed form solution of the following equation
\begin{align}
y(x)-y(x+1)=
\frac{2x+a}{(x^2+a x+b_1)(x^2+a x+b_2)}.
\end{align}

\subsection{The rational series $\sum\frac{2(p n + q)}{(p n + q)^4 + s (p n + q)^2 + r}$}

Firstly, we shall study  the following difference equation of order one
\begin{align}
y(x)-y(x+1)=\frac{2(p x + q)}{(p x + q)^4 + s (p x + q)^2 + r}.
\end{align}
\begin{thm}
Let $p, q, r, s\in\mathbb{R}$ with $p>0$. We define $H_2(p,q,r,s;x)$~(or shortly $H_2(x)$) by
\begin{align}
H_2(p,q,r,s;x):=\frac{\frac{1}{ p^3}}{(x+\omega)^2+\frac 14+\frac{ s}{2p^2}+\K_{n=1}^{\infty}\left(\frac{\kappa_n}
{(x+\omega)^2+\lambda_n}\right)},
\end{align}
where
\begin{align*}
\omega = -\frac 12+\frac{q}{ p},\quad \kappa_n=\frac{n^2 \left(- n^4  - \frac{2s}{p^2} n^2 +\frac{ 4 r- s^2}{p^4}\right)}{
2^2 (2 n - 1) (2 n + 1)},\quad\lambda_n=\frac{2 n^2 + 2 n + 1}{4} + \frac{s}{2 p^2}.
\end{align*}
Let $x\notin\{l+\alpha: l\in\mathbb{Z}\backslash\mathbb{N},~ (p\alpha+q)^4+s(p\alpha+q)^2 +r=0, ~\alpha\in\mathbb{C} \}$. If either one of $\sqrt{-2 \sqrt{r} - s}/p$, $\sqrt{2 \sqrt{r} - s}/p$ is a positive integer, or $\Re x>-\omega$, then
\begin{align}
H_2(x)-H_2(x+1)=\frac{2(p x + q)}{(p x + q)^4 + s (p x + q)^2 + r}.\label{Difference Equation-6}
\end{align}
\end{thm}

\bigskip

\proof Firstly, we note that
\begin{align}
n^4 +\frac{2s}{p^2} n^2 +\frac{s^2 -4 r }{p^4}
=\left(\frac{-s+2\sqrt{r}}{p^2}-n^2\right)
\left(\frac{-s-2\sqrt{r}}{p^2}-n^2\right).
\end{align}
We shall discuss seven cases.
\bigskip

\noindent\emph{(Case 1)} $r\ge 0$ and $ -2 \sqrt{r} - s> 0$.
In this case, we have $ 2 \sqrt{r} - s> 0$.
Applying Lemma 2 with $(x,l,m)=(2(x+\omega),\sqrt{2 \sqrt{r} - s}/p,\sqrt{-2 \sqrt{r} - s}/p)$ and dividing both sides by $2n \sqrt{s^2-4r}/p^2$, we assume that the conditions of Theorem 9 hold, then
\begin{align}
\frac{p^2}{2 \sqrt{s^2-4r}}\frac{1-P}{n(1+P)}
=&\frac{1}{\left(2(x+\omega)\right)^2-n^2+2s/p^2+1
+\K_{j=1}^{\infty}\left(
\frac{4\left((2 \sqrt{r} - s)/p^2-j^2\right)\left((-2 \sqrt{r} - s)/p^2-j^2\right)(n^2-j^2)}{(2j+1)
\left(\left(2(x+\omega)\right)^2-n^2+2s/p^2
+2j^2+2j+1\right)}\right)}\\
=&\frac{1}{\left(2(x+\omega)\right)^2-n^2+2s/p^2+1
+\K_{j=1}^{\infty}\left(
\frac{4\left(j^4+2sj^2/p^2+(s^2-4r)/p^4\right)(n^2-j^2)}{(2j+1)
\left(\left(2(x+\omega)\right)^2-n^2+2s/p^2
+2j^2+2j+1\right)}\right)}.\nonumber
\end{align}
Now let $n$ tend to zero. On the right side, we arrive at
\begin{align}
\frac{1}{\left(2(x+\omega)\right)^2+2s/p^2+1
+\K_{j=1}^{\infty}\left(
\frac{-4\left(j^4+2sj^2/p^2+(s^2-4r)/p^4\right)j^2}{(2j+1)
\left(\left(2(x+\omega)\right)^2+2s/p^2
+2j^2+2j+1\right)}\right)}=\frac {p^3}{4} H_2(x).
\end{align}
On the other hand, from the definition of $P$, it is easy to see that $\lim_{n\rightarrow 0}P=1$. A direct calculation with the use of L'Hospital's rule gives
\begin{align}
&\lim_{n\rightarrow 0}\frac{1-P}{n(1+P)}=\lim_{n\rightarrow 0}
\frac{1}{1+P}\lim_{n\rightarrow 0}\frac{1-P}{n}
\lim_{n\rightarrow 0}=\frac 12\lim_{n\rightarrow 0}
\frac {\partial}{\partial n}(1-P)\\
=&\frac{1}{2 }\left\{-\psi\left(x+\frac qp+\frac {1}{2p}(-\sqrt{-2 \sqrt{r} - s}-\sqrt{2 \sqrt{r} - s})\right)+\psi\left(x+\frac qp+\frac{1}{2p}(\sqrt{-2 \sqrt{r} - s}-\sqrt{2 \sqrt{r} - s})\right)\right.\nonumber\\
&\left.+\psi\left(x+\frac qp+\frac{1}{2p}(-\sqrt{-2 \sqrt{r} - s}+\sqrt{2 \sqrt{r} - s})\right)-\psi\left(x+\frac qp+\frac{1}{2p}(\sqrt{-2 \sqrt{r} - s}+\sqrt{2 \sqrt{r} - s})\right)\right\}\nonumber\\
=&\frac 12\sum_{k=0}^{\infty}\left(\frac{1}{k+x+\frac qp+\frac {1}{2p}(-\sqrt{-2 \sqrt{r} - s}-\sqrt{2 \sqrt{r} - s})}-\frac{1}{k+x+\frac qp+\frac {1}{2p}(\sqrt{-2 \sqrt{r} - s}-\sqrt{2 \sqrt{r} - s})}\right.\nonumber\\
&\left.\quad -\frac{1}{k+x+\frac qp+\frac{1}{2p} (-\sqrt{-2 \sqrt{r} - s}+\sqrt{2 \sqrt{r} - s})}+\frac{1}{k+x+\frac qp+\frac{1}{2p}(\sqrt{-2 \sqrt{r} - s}+\sqrt{2 \sqrt{r} - s})}\right)\nonumber\\
=&p\sqrt{s^2-4r}\sum_{k=0}^{\infty}\frac{\left(p(k+x)+q\right)}
{\left(p(k+x)+q\right)^4+s\left(p(k+x)+q\right)^2+r}.
\nonumber
\end{align}
Hence,
\begin{align}
&H_2(x)-H_2(x+1)\\
=&\sum_{k=0}^{\infty}\frac{2\left(p(k+x)+q\right)}
{\left(p(k+x)+q\right)^4+s\left(p(k+x)+q\right)^2+r}
-\sum_{k=0}^{\infty}\frac{2\left(p(k+1+x)+q\right)}
{\left(p(k+1+x)+q\right)^4+s\left(p(k+1+x)+q\right)^2+r}\nonumber\\
=&\frac{2\left(px+q\right)}
{\left(px+q\right)^4+s\left(px+q\right)^2+r}.\nonumber
\end{align}
This proves~\eqref{Difference Equation-6} in Case 1.
\bigskip

\noindent\emph{(Case 2)} $r\ge 0$, $2 \sqrt{r} - s> 0$, and
$ -2 \sqrt{r} - s<0$. Applying Lemma 2 with $(x,l,m)=(2(x+\omega),\sqrt{2 \sqrt{r} - s}/p,\sqrt{2 \sqrt{r} + s}~i/p)$, and then let $n$ tend to zero. The proof is very similar to that of Case 1, we omit the detail here.
\bigskip

\noindent\emph{(Case 3)} $r\ge 0$ and $ 2 \sqrt{r} - s< 0$. In this case, we have $ -2 \sqrt{r} - s<0$. Applying Lemma 2 with $(x,l,m)=(2(x+\omega),\sqrt{-2 \sqrt{r} + s}~i/p,\sqrt{2 \sqrt{r} + s}~i/p)$, and let $n$ tend to zero. The proof is very similar to that of Case 1.

\bigskip

\noindent\emph{(Case 4)} $r> 0$ and $ s=2 \sqrt{r}$. Applying Lemma 2 with $(x,l)=(2(x+\omega),\sqrt{4 \sqrt{r}}~i/p)$ and dividing both sides by $2m n \sqrt{4 \sqrt{r}}~i/p$, we have
\begin{align}
&\frac{p}{4 r^{1/4}~i}\frac{1-P}{m n(1+P)}\\
=&\frac{1}{\left(2(x+\omega)\right)^2-m^2-n^2+4 \sqrt{r}/p^2+1
+\K_{j=1}^{\infty}\left(
\frac{-4\left(4 \sqrt{r}/p^2+j^2\right)(m^2-j^2)(n^2-j^2)}{(2j+1)
\left(\left(2(x+\omega)\right)^2-m^2-n^2+4 \sqrt{r}/p^2
+2j^2+2j+1\right)}\right)}.\nonumber
\end{align}
We let $n$ and $m$ tend to zero, successively. On the right side, one has
\begin{align}
\frac{1}{\left(2(x+\omega)\right)^2+4 \sqrt{r}/p^2+1
+\K_{j=1}^{\infty}\left(
\frac{-4\left(4 \sqrt{r}/p^2+j^2\right)j^4}{(2j+1)
\left(\left(2(x+\omega)\right)^2+4 \sqrt{r}/p^2
+2j^2+2j+1\right)}\right)}
=\frac{p^3}{4}H_2(x).
\end{align}
By using of L'Hospital's rule, we deduce that
\begin{align}
&\lim_{n\rightarrow 0}\frac{1-P}{n(1+P)}=\lim_{n\rightarrow 0}
\frac{1}{1+P}\lim_{n\rightarrow 0}\frac{1-P}{n}
\lim_{n\rightarrow 0}=\frac 12\lim_{n\rightarrow 0}
\frac {\partial}{\partial n}(1-P)\\
=&\frac{1}{2 }\left\{-\psi\left(-\frac m2+\frac{q-r^{1/4}~i}{p}+x\right)
+\psi\left(\frac m2+\frac{q-r^{1/4}~i}{p}+x\right)\right.\nonumber\\
&\left.\quad+\psi\left(-\frac m2+\frac{q+r^{1/4}~i}{p}+x\right)
-\psi\left(\frac m2+\frac{q+r^{1/4}~i}{p}+x\right)
\right\}.\nonumber
\end{align}
By making use of L'Hospital's rule again, we find that
\begin{align}
&\lim_{m\rightarrow 0}\lim_{n\rightarrow 0}\frac{1-P}{mn(1+P)}
=\lim_{m\rightarrow 0}\frac 1m \left(\lim_{n\rightarrow 0}\frac{1-P}{n(1+P)}\right)\\
=&\frac 12\left\{\psi'\left(x+\frac{q-r^{1/4}~i}{p}\right)
-\psi'\left(x+\frac{q+r^{1/4}~i}{p}\right)\right\}\nonumber\\
=&\frac 12\sum_{k=0}^{\infty}\frac{1}{\left(k+x+\frac{q-r^{1/4}~i}{p}
\right)^2}
-\frac 12\sum_{k=0}^{\infty}\frac{1}{\left(k+x+\frac{q+r^{1/4}~i}{p}
\right)^2}\nonumber\\
=&2p^2 r^{1/4}~i\sum_{k=0}^{\infty}\frac{\left(p(k+x)+q\right)}
{\left(\left(p(k+x)+q\right)^2
+\sqrt{r}\right)^2}\nonumber\\
=&p^2 r^{1/4}~i\sum_{k=0}^{\infty}
\frac{2\left(p(k+x)+q\right)}{\left(p(k+x)+q\right)^4
+s\left(p(k+x)+q\right)^2
+r}.\nonumber
\end{align}
Following the same argument as Case 1, we find that Theorem 9 holds in this case.
\bigskip

\noindent\emph{(Case 5)} $r> 0$ and $ s=-2 \sqrt{r}$. Applying Lemma 2 with $(x,l)=(2(x+\omega),\sqrt{4 \sqrt{r}}/p)$, and let $n$ and $m$ tend to zero, successively. The proof is very similar to that of Case 4.
\bigskip

\noindent\emph{(Case 6)} $r<0$. Applying Lemma 2 with $(x,l,m)=(2(x+\omega),\sqrt{2 \sqrt{-r}~i - s}/p,\sqrt{-2 \sqrt{-r}~i - s}/p)$, and let $n$ tend to zero. The proof is very similar to that of Case 1.
\bigskip

\noindent\emph{(Case 7)} $r=s=0$. By the case 3 in the proof of Theorem 3 with $a=3q/p$, we deduce that Theorem 9 holds true.

Finally, combining Case 1 to 7 will finish the proof of Theorem 9.\qed

\bigskip

\begin{thm} With the conditions of Theorem 9, let $n_2$ be a non-negative integer such that $n_2>\max\limits_{\alpha}\{\alpha:~ (p\alpha+q)^4+s(p\alpha+q)^2 +r=0,~\alpha\in\mathbb{Z}\}$, and either one of $\sqrt{-2 \sqrt{r} - s}/p$, $\sqrt{2 \sqrt{r} - s}/p$ is a positive integer, or $n_2>\omega$, then
\begin{align}
\sum_{n=n_2}^{\infty}\frac{2(pn+q)}{(pn+q)^4+s(p n+q)^2+r}=H_2(p,q,s,r;n_2).\label{AN extended MS-2}
\end{align}
In particular, let
\begin{align}
r=\frac{(p^2k^2 + s)^2}{4},\quad k\in\mathbb{N}.\label{AN extended MS-2-1}
\end{align}
Assume that $n_2>\max\limits_{\alpha}\{-1, \alpha:~ (p\alpha+q)^4+s(p\alpha+q)^2 +r=0,~\alpha\in\mathbb{Z}\}$, then
\begin{align}
\sum_{n=n_2}^{\infty}\frac{2(pn+q)}{(pn+q)^4+s (pn+q)^2+r}
=
\frac{1}{\left(n_2-\frac 12+\frac qp\right)^2+\frac{1 }{4}+\frac{s}{2p^2}+\K_{n=1}^{k-1}\left(
\frac{\frac{n^2 (n^2+k^2+2s/{p^2})(k^2-n^2)}{4 (2 n - 1)
(2 n + 1)}}{\left(n_2-\frac 12+\frac qp\right)^2+\frac{2 n^2+ 2 n +1}
{4}+\frac{s}{2p^2}}\right)}.
\end{align}
\end{thm}
\proof It follows from Theorem 9 and the telescoping method readily.\qed

\noindent\textbf{Example 8} Taking $(p,q,s,r)=(1,0,2r^2,r^4)$ in Theorem 10~(or $(a,b_1,b_2)=(0,r^2,r^2)$ in \eqref{AN extended MS-1}), the series in \eqref{AN extended MS-2} become the Mathieu series. Hence, for $l\ge 1$
\begin{align}
S(r)=\sum_{m=1}^{\infty}\frac{2m}{(m^2+r^2)^2}
=\sum_{m=1}^{l-1}\frac{2m}{(m^2+r^2)^2}+H_2(1,0,2r^2,r^4;l).
\end{align}

\bigskip

\noindent\textbf{Example 9} Let $\lambda$ be real. Consider the infinite series
\begin{align}
T(\lambda)=\sum_{n=1}^{\infty}\frac{2n}{n^4+\lambda n^2+\frac{ (4+\lambda)^2}{4}}.
\end{align}
We take $k=2$ and $(p,q,s)=(1,0,\lambda)$ in~\eqref{AN extended MS-2-1}, then
\begin{align}
H_2\left(1,0,\lambda, \frac{ (4+\lambda)^2}{4};x\right)=\frac{2 (3 - 2 x + 2 x^2 + \lambda)}{(2 + 2 x^2 + \lambda) (4 - 4 x + 2 x^2 +\lambda)},
\end{align}
and for $\lambda\neq -2, -4$
\begin{align}
T(\lambda)=H_2\left(1,0,\lambda, \frac{ (4+\lambda)^2}{4};1\right)=\frac{2 (3 + \lambda)}{(2 + \lambda) (4 + \lambda)}.
\end{align}
Note that
\begin{align}
T(\sqrt{2})=\frac{2 (3 + \sqrt{2})}{10 + 6 \sqrt{2}}=\frac{9 - 4 \sqrt{2}}{7},
\end{align}
hence $T(\sqrt{2})$ is an irrational number. In fact, the assertion above may be checked directly by hands.


\section{Two applications of Theorem 8 and 10}

\subsection{The alternating Mathieu series}
\begin{thm} Let the formal continued fraction $H_1(a,b_1,b_2;x)$ be defined by~\eqref{H1-Def}. For all positive integer $k_1$ and $k_2$, we have
\begin{align}
\tilde{S}(r)=&\frac 18\sum_{m=0}^{k_1-1}\frac{2m+1}{\left(m^2+m+\frac {1+r^2}{4}\right)^2}-\frac 18\sum_{m=0}^{k_2-1}\frac{2m}{\left(m^2+(\frac r2)^2\right)^2}
\label{A-Mathiew-CF}\\
&+
\frac 18 H_1\left(1,\frac {1+r^2}{4},\frac {1+r^2}{4};k_1\right)-\frac 18 H_1\left(0,\frac {r^2}{4},\frac {r^2}{4};k_2\right).\nonumber
\end{align}
In particular,
\begin{align}
\tilde{S}(r)=\frac{2}{(1+r^2)^2}+\frac 18 H_1\left(1,\frac {1+r^2}{4},\frac {1+r^2}{4};1\right)-\frac 18 H_1\left(0,\frac {r^2}{4},\frac {r^2}{4};1\right).\label{A-Mathiew-CF-1}
\end{align}
\end{thm}
\proof From the definition of the alternating Mathieu series in~\eqref{Alternating Mathieu-Series}, we rewrite it into two parts
\begin{align}
\tilde{S}(r)=&\sum_{m=1}^{\infty}(-1)^{m-1}\frac{2m}{(m^2+r^2)^2}
=\sum_{n=0}^{\infty}\frac{2\cdot (2n+1)}{\left((2n+1)^2+r^2\right)^2}-\sum_{n=1}^{\infty}
\frac{2\cdot 2n}{\left((2n)^2+r^2\right)^2}\\
=&\frac 18\sum_{n=0}^{\infty}\frac{2n+1}{\left(n^2+n+\frac {1+r^2}{4}\right)^2}-\frac 18\sum_{n=0}^{\infty}\frac{2n}{\left(n^2+(\frac r2)^2\right)^2}.\nonumber
\end{align}
Applying Theorem 8 with $(a,b_1,b_2)=(1,\frac {1+r^2}{4},\frac {1+r^2}{4})$ and $(a,b_1,b_2)=(0,\frac {r^2}{4},\frac {r^2}{4})$, respectively,  we get the desired assertion. Finally, on taking $k_1=k_2=1$, \eqref{A-Mathiew-CF-1} follows from \eqref{A-Mathiew-CF} readily.\qed

\subsection{For the numbers $\mathbb{M}_2^{(m,j)}$}
Let $m\in\mathbb{N}$ and $j\in\{1,2,\ldots,m-1\}$ with $(m,j)=1$. P.~J. Szablowski~\cite{Sza} introduced the numbers
\begin{align}
\mathbb{M}_k^{(m,j)}=\sum_{n=0}^{\infty}\frac{(-1)^n}{(mn+j)^k}.
\end{align}
Notice that $\mathbb{M}_k^{(1,1)}=\sum_{j=1}^{\infty}
(-1)^{j-1}/j^k$ and $\mathbb{M}_1^{(2,1)}=\pi/4$. The number $\mathbb{M}_2^{(2,1)}=\sum_{j=0}^{\infty}\frac{(-1)^{j}}{(2j+1)^2}$
is Catalan constant $\mathbb{K}$, which is one of those classical constants whose irrationality and transcendence remain unproven. Two continued fraction formulas for Catalan constant can be found in Berndt~\cite[p.~151\nobreakdash--153]{Ber}. W. Zudilin~\cite{Zud} obtained new one, also see Cuyt et al.~\cite[Eq.~(10.12.5), p.~189]{CPV}. When $k=1$, the continued fraction representation of the numbers $\mathbb{M}_1^{(m,j)}$ could be obtained from Theorem 2 easily. When $k=3$, quite similarly to the alternating Mathieu series in Theorem 11, we may use Theorem 4 to write $\mathbb{M}_3^{(m,j)}$ into a linear combination of two continued fractions. Now we state the main result as follows.

\begin{thm} Let $\omega=\frac{2 j - 3m}{4 m}$, and the formal continued fraction $CF_2(m,j;x)$ be defined by
\begin{align}
CF_2(m,j;x):=\frac{1}{(x+\omega)^2+\frac{3}{16}+
\K_{n=1}^{\infty}\left(\frac{-\frac{n^4}{16}}{(x+\omega)^2
+\frac{8 n^2 + 8 n + 3}{16}}\right)}.
\end{align}
We let $a=\frac jm - \frac 12, b=-\frac{j}{4 m} + \frac{j^2}{4 m^2}$. For all positive integer $l$, then
\begin{align}
\mathbb{M}_2^{(m,j)}=\frac{1}{j^2}-\frac{1}{8m^2}
\sum_{n=1}^{l-1}\frac{2n+b}{(n^2+an+b)^2}
-\frac{1}{8m^2}CF_2(m,j;l).
\end{align}
In particular,
\begin{align}
\mathbb{M}_2^{(m,j)}=\frac{1}{j^2}-\frac{1}{8m^2}\frac{1}
{\left(\frac{2j+m}{4m}\right)^2+\frac{3}{16}+
\K_{n=1}^{\infty}\left(\frac{-\frac{n^4}{16}}
{\left(\frac{2j+m}{4m}\right)^2+\frac{8 n^2 + 8 n + 3}{16}}\right)}.
\end{align}
\end{thm}

\proof First, we note that the equalities $0<\frac{3m-2j}{4m}<1$ always hold. Following the same argument as Theorem 11, we also have
\begin{align}
\mathbb{M}_2^{(m,j)}=\frac{1}{j^2}-\frac{1}{8m^2}
\sum_{n=1}^{\infty}\frac{2n+b}{(n^2+an+b)^2}.
\end{align}
It is elementary to check that
\begin{align*}
\kappa_n=\frac{n^4 (a^2 - 4 b - n^2)}{4 (2 n - 1) (2 n + 1)}=-\frac{n^4}{16},\quad
\lambda_n=\frac{2 n^2 + 2 n + 1 + 4 b - a^2}{4}=\frac{8 n^2 +
8 n + 3}{16}.
\end{align*}
Now Theorem 12 follows easily from Theorem 8 .\qed

\section{Conclusions}
From the above discussion, we observe that for a specific rational series, the \emph{multiple-correction method} provides a useful tool for finding a simple closed form solution, testing and guessing the continued fraction representation, or getting the fastest possible finite continued fraction approximation, etc. So our method should help advance the approximation theory, the theory of continued fraction and the generalized hypergeometric function, etc. Furthermore, probably these continued fraction formulas could be used to study the irrationality, transcendence of the involved series.

\begin{flushleft}
\end{flushleft}

\end{document}